\documentclass[11pt,a4paper]{article}
\usepackage{amsmath}
\usepackage{amsfonts}
\usepackage{amssymb}
\usepackage{amsthm}
\usepackage{mathrsfs}
\usepackage{dsfont}
\usepackage{paralist}
\usepackage{natbib}
\usepackage{bm}
\usepackage{color}
\usepackage[colorlinks=true,linkcolor=blue,citecolor=blue,pdfborder={0 0 0}]{hyperref}
\usepackage{graphicx}
\usepackage{tabularx}
\usepackage[left=3.5cm,right=3.5cm,top=2cm,bottom=2cm,includeheadfoot]{geometry}
\usepackage{comment}
\usepackage{tikz}
\usetikzlibrary{decorations.pathreplacing}
\usetikzlibrary{decorations.pathmorphing}
\usepackage{resizegather}
\DeclareFontFamily{OT1}{pzc}{}
\DeclareFontShape{OT1}{pzc}{m}{it}{<-> s * [1.10] pzcmi7t}{}
\DeclareMathAlphabet{\mathpzc}{OT1}{pzc}{m}{it}

\newtheoremstyle{normal}
{2ex}               
{3ex}               
{}                  
{}                  
{\bfseries} 
{}                  
{2pt}   
{\thmname{#1}\thmnumber{ #2.} \thmnote{(#3)}}

\newtheoremstyle{italic}
{2ex}
{3ex}
{\itshape}
{}
{\bfseries} 
{}
{2pt}
{\thmname{#1}\thmnumber{ #2.} \thmnote{(#3)}}

\theoremstyle{normal}
\newtheorem{definition}{Definition}[section]
\newtheorem{remark}[definition]{Remark}
\newtheorem{example}[definition]{Example}
\newtheorem{condition}[definition]{Condition}

\theoremstyle{italic}
\newtheorem{theorem}[definition]{Theorem}
\newtheorem{lemma}[definition]{Lemma}

\renewcommand{\P}{\mathbb{P}}

\newcommand{\E}{\mathbb{E}}

\DeclareMathAccent{\verywidehat}{\mathord}{largesymbols}{'144}

\newcommand{\Q}{\mathbb{Q}}
\newcommand{\N}{\mathbb{N}}
\newcommand{\R}{\mathbb{R}}

\newcommand{\al}{\alpha}

\newcommand{\la}{\lambda}

\newcommand{\ga}{\gamma}
\newcommand{\ka}{\kappa}

\newcommand{\si}{\sigma}

\newcommand{\eps}{\varepsilon}
\newcommand{\vpi}{\varpi}

\newcommand{\de}{\delta}
\newcommand{\te}{\theta}

\newcommand{\De}{\Delta}
\newcommand{\om}{\omega}
\newcommand{\Om}{\Omega}

\newcommand{\fF}{\mathcal F}

\newcommand{\lL}{\mathcal L}

\newcommand{\proba}{(\Omega ,\fF,(\fF_t)_{t\geq0},\P)}

\newcommand{\toop}{\stackrel{\P}{\longrightarrow}}
\newcommand{\tols}{~\stackrel{\lL-(s)}{\longrightarrow}~}
\newcommand{\tol}{\stackrel{\lL}{\longrightarrow}}

\newcommand{\pn}{\stackrel{\P}{\longrightarrow}}
\newcommand{\weak}{\stackrel{w}{\longrightarrow}}

\newcommand{\wte}{\widehat{\theta}}
\newcommand{\wSi}{\widehat{\Sigma}}

\newcommand{\wV}{\widehat{V}}

\newcommand{\wY}{\widehat{Y}}

\newcommand{\WV}{\widetilde{V}}

\newcommand{\Wsi}{\widetilde{\sigma}}

\newcommand{\Wb}{\widetilde{b}}

\newcommand{\Ind}{\mathds{1}}

\DeclareMathOperator{\Var}{Var}
\DeclareMathOperator{\Cov}{Cov}

\allowdisplaybreaks[1]

\begin{document}

\title{A universal approach to estimate the conditional variance in semimartingale limit theorems}

\author{Mathias Vetter\thanks{Christian-Albrechts-Universit\"at zu Kiel, Mathematisches Seminar, Ludewig-Meyn-Str.\ 4, 24118 Kiel, Germany.
{E-mail:} vetter@math.uni-kiel.de.} \bigskip \\
{Christian-Albrechts-Universit\"at zu Kiel}
}

\maketitle

\begin{abstract}
The typical central limit theorems in high-frequency asymptotics for semimartingales are results on stable convergence to a mixed normal limit with an unknown conditional variance. Estimating this conditional variance usually is a hard task, in particular when the underlying process contains jumps. For this reason, several authors have recently discussed methods to automatically estimate the conditional variance, i.e.\ they build a consistent estimator from the original statistics, but computed at various different time scales. Their methods work in several situations, but are essentially restricted to the case of continuous paths always. The aim of this work is to present a new method to consistently estimate the conditional variance which works regardless of whether the underlying process is continuous or has jumps. We will discuss the case of power variations in detail and give insight to the heuristics behind the approach. 
\end{abstract}

\medskip

\textit{Keywords and Phrases:} Asymptotic conditional variance; high-frequency statistics; It\^o semimartingale; jumps; stable convergence

%


\section{Introduction}
\def\theequation{1.\arabic{equation}}
\setcounter{equation}{0}

The asymptotic theory for functionals of semimartingales observed at high frequency is well understood now. Since the beginning of the century a variety of laws of large numbers and accompanying central limit theorems has been stated in different situations, starting with power and bipower variation of continuous processes (\cite{barshe2003} or \cite{baretal2006}). Crucial generalizations involve the case of possible jumps in the process (\cite{jacod2008}) or the discussion of observations with additional microstructure noise (\cite{jacodetal2010}). Later extensions regard truncated increments, multivariate processes or the treatment of irregularity and asynchronicity in the data. A general overview about these results and statistical applications can be found in the monographs \cite{JacPro12} and \cite{AitJac14}. 

Typically the central limit theorems in these situations are stated as follows: One proves stable convergence in law of an appropriately rescaled statistic to a mixed normal limit, where the (asymptotic) conditional variance of the limiting variable is a random variable which depends in a complicated way on the underlying semimartingale. Once a consistent estimator for this conditional variance has been constructed, thanks to the properties of stable convergence in law, one can deduce the convergence in distribution of the standardized statistic to a standard Gaussian law. This opens the door for all kinds of statistical applications. 

Constructing a consistent estimator for the conditional variance, however, is not always a simple task. Compared with the original object of interest for which the law of large numbers is shown, usually an integral of a power of volatility or a sum of a power of jumps, the variance is typically of a more complicated form and might depend on additional objects as well. In particular, apart from the case of power variations of continuous processes, it is not possible to estimate the variance by using similar statistics as for the corresponding law of large numbers. Hence, estimators are usually constructed based on the specific form of the conditional variance in the respective situations. This procedure has two major drawbacks: First, every newly proven central limit theorem requires new estimators for the conditional variances. Second, when the model is not correctly specified, it is likely that the proposed estimator does not work. 

A different approach is to build an estimator which only requires knowledge of the original statistics and does not rely on the specific form of the conditional variance. For example, \cite{jacod2008} discusses statistics of the form
\[
U_n = \sum_{i=1}^n f_n(\De_i^n X),
\]
for simplicity over $[0,1]$, where $\De_i^n X = X_{i \De_n} - X_{(i-1)\De_n}$ denotes the $i$-th increment of the semimartingale $X$, $\De_n \to 0$, and where $f_n: \R \to \R$ is a function which may or may not depend on $n$. Several laws of large numbers and associated central limit theorems are shown in various cases. A universal estimator for the conditional variance in these very central limit theorems would then only depend on $f_n$, but not utilize the specific form of the conditional variance in the respective situations. Whether such estimators exist, and how they look like, is obviously an important question in the theoretical discussion of high-frequency statistics. 

In recent years two classes of such universal estimators have been proposed in the literature. \cite{mykzha2017} base their estimator on a comparison of local versions of $U_n$ computed over neighbouring intervals of length $k_n \De_n$, $k_n \to \infty$ and $k_n \De_n \to 0$, whereas \cite{chrietal2017} use a subsampling approach which compares $U_n$ with versions where only every $k_n$-th increment is taken into account. Both estimators are shown to work in a variety of situations, but only when the semimartingale $X$ does not jump (or when the jumps do not contribute to the limiting distribution), and it is rather simple to see that both procedures indeed do not work when the limiting distribution contains jumps. 

Therefore, the question remains whether it is possible to construct a universal estimator for the conditional variance which works both in the continuous case and in the case involving jumps and, if yes, how it could be constructed. We will give positive answers to both questions, for simplicity in the case of power variations only, which means that $X$ is a general Ito semimartingale including jumps and that $f_n$ is essentially of the form $f_n(x) = |x|^p$, $p > 0$, up to a possible standardization. Already in this situation we will see all different kinds of limiting behaviour, including conditional variances which only depend on the volatility or which depend jointly on jumps and volatility. It is to be expected that the same construction of a universal estimator works for most other statistics as well, as the main idea behind the proof of the respective central limit theorems usually is the same as for the corresponding power variations. 
 
The paper is organized as follows: After introducing the setting in Section \ref{sec:setting}, we will discuss three novel universal estimators for the conditional variance in Section \ref{sec:results}. While the first two estimators are rather simple to construct in practice, they have the deficiency that they do not work in all situations. In fact, the first one is consistent for continuous processes, but when jumps dominate it only converges stably in law to a random variable whose mean is the conditional variance. Similarly for the second estimator, but with different roles. The estimator is consistent in the jump case, but does not converge to the correct conditional variance for continuous processes. A remarkable exception is the case $p=2$ in which it gives an alternative estimator for the conditional variance when the quadratic variation is to be estimated. Finally, the intuition behind both estimators is combined to construct the universal estimator which formally works in all situations. Its computation time is of order $n {\binom{k_n}{\ell_n}}$ for sequences $k_n$ and $\ell_n$ converging to infinity, however, so it is of theoretical interest in the first place rather than being a serious alternative in all practical cases. The proofs are given in Section \ref{sec:proof}.  

\section{Setting} \label{sec:setting}
\def\theequation{2.\arabic{equation}}
\setcounter{equation}{0}

Suppose that we have a filtered probability space $\proba$ on which an Ito semimartingale of the form 
\begin{multline}
\label{ItoSemimart}
X_t = X_0 + \int \limits_0^t b_s ds + \int \limits_0^t \sigma_s dW_s + \int \limits_0^t \int \delta(s,z)\mathds{1}_{\{|\delta(s,z)|\leq 1\}} (\mu - \nu)(ds,dz) \\
+ \int \limits_0^t \int \delta(s,z) \mathds{1}_{\{|\delta(s,z)|> 1\}} \mu(ds,dz)
\end{multline}
is defined, where $W$ is a standard Brownian motion, $\mu$ is a Poisson random measure on $\mathbb{R}^+ \times\mathbb{R}$, and its predictable compensator satisfies \mbox{$\nu(ds,dz)=ds \otimes \lambda(dz)$} for some $\sigma$-finite measure $\lambda$ on $\mathbb{R}$ endowed with the Borelian $\sigma$-algebra. We further assume that $b$ and $\sigma$ are adapted processes and that $\delta$ is predictable on $\Omega \times \mathbb{R}^+ \times \mathbb{R}$. We write $\Delta X_s=X_s-X_{s-}$ with $X_{s-}=\lim_{t \nearrow s} X_t$ for a possible jump of $X$ in $s$. 

We will work in a high-frequency framework, so without loss of generality we assume to be on the fixed interval $[0,1]$. Observations of $X$ take place at the regular times $i \De_n$, $i=0, \ldots, n$, where we set $n = \De_n^{-1}$. Throughout the paper, $\De_n \to 0$ governs the asymptotics. 

In order to prove asymptotic results for statistics based on increments of $X$, one typically needs additional assumptions on the semimartingale characteristics. Our aim in the following is not to be as general as possible, so we will state sufficient conditions in order to prove consistency of the statistics and associated central limit theorems, respectively. The first one is good enough for theorems on consistency, and it even is sufficient for some central limit theorems. 

\begin{condition} \label{condlln}
The process $(b_s)$ is locally bounded and predictable, the process $(\si_s)$ is c\`adl\`ag, and  there exist a sequence $(\tau_n)$ of stopping times increasing to infinity and a sequence $(\ga_n)$ of deterministic real functions such that $1 \wedge |\delta(s,z)| \le \ga_n(z)$ for all $s \le \tau_n$ and $\int \ga_n(z)^2 \la(dz) < \infty$ hold. 
\end{condition}

Stronger assumptions are typically needed when one is interested in a central limit theorem accociated with a limit in probability which is governed by the continuous martingale part of $X$. What is always needed is that $\si$ is positive and that it takes a form similar to (\ref{ItoSemimart}). 

\begin{condition} \label{condclt}
We assume that the process $(\si_s)$ is bounded below by a positive number and of the form 
\[
\si_t = \si_0 + \int_0^t \Wb_s ds +  \int_0^t \Wsi_s dW_s + M_t + \sum_{0 < s \le t} \De\si_s \Ind_{\{|\De \si_s| > 1\}}
\] 
with $M$ being a local martingale with $|\De M_s| \le 1$, orthogonal to $W$, and we assume that $\langle M, M \rangle_t = \int_0^t \al_s ds$ as well as that the compensator of $\sum_{0 < s \le t} \De\si_s \Ind_{\{|\De \si_s| > 1\}}$ takes the form $\int_0^t \al'_s ds$. The processes $(b_s)$ and $(\Wsi_s)$ are c\`adl\`ag, and the processes $(\Wb_s)$, $(\al_s)$ and $(\al'_s)$ are locally bounded and predictable.  
\end{condition}

Even this condition is not general enough in the case where $X$ has jumps as well; see Theorems 5.3.5 and 5.3.6 in \cite{JacPro12}. We will therefore assume that $X$ is continuous whenever we are concerned with central limit theorems associated to the continuous martingale part only. Condition \ref{condclt} turns out to be sufficient then. 

\section{Results} \label{sec:results}
\def\theequation{3.\arabic{equation}}
\setcounter{equation}{0}

\subsection{Limit theorems for power variations}
The typical object of interest in high-frequency statistics is a statistic of the form 
\[
U_n = \sum_{i=1}^n f_n(\De_i^n X),
\]
where $\De_i^n X = X_{i \De_n} - X_{(i-1)\De_n}$ denotes the $i$-th increment of $X$ and $f_n: \R \to \R$ is a function which may or may not depend on $n$. Typical examples are power variations of the form
\[
f_n(x) = |x|^p \quad \text{or} \quad f_n(x) = \Delta_n^{1-p/2} |x|^p
\]
for some $p > 0$, where the latter scaling depends on the length of the interval over which the increment is computed. For those power variations and related statistics, a rule of thumb is: Whenever a weak law of large numbers holds, the limit is of the form
\[
U = \int_0^1 g(\si_s) ds + \sum_{0 < s \le 1} h(\De X_s),
\]
where $g: \R^+ \to \R$ and $h: \R \to \R$ are suitable functions depending on $f_n$. Let us recall the results from Theorem 2.2 and Theorem 2.4 in \cite{jacod2008}. 

\begin{theorem} \label{jacodlln}
Let $X$ be a semimartingale of the form (\ref{ItoSemimart}) and assume that Condition \ref{condlln} holds.
\begin{itemize}
	\item[(a)] Let $p < 2$ and $f_n(x) = \Delta_n^{1-p/2} |x|^p$. Then
	\[
	U_n \toop m_p \int_0^1 \si_s^p ds
	\]
	with $m_p = \E[|N|^p]$ for $N \sim \mathcal N(0,1)$. 
	\item[(b)] Let $p > 2$ and $f_n(x) = |x|^p$ for any $n$. Then 
	\[
	U_n \toop \sum_{0 < s \le 1} |\De X_s|^p.
	\]
	\item[(c)] Let $f_n(x) = |x|^2$ for any $n$. Then 
	\[
	U_n \toop [X,X]_1 = \int_0^1 \si_s^2 ds + \sum_{0 < s \le 1} |\De X_s|^2.
	\] 
\end{itemize}
\end{theorem}

\begin{remark} \label{remlln}
In the case where no jumps are present, the law of large numbers in part (a) also holds for $p \ge 2$. Similarly, if the continuous martingale part vanishes the claim in part (b) also holds for $p \in (1,2]$ and, under a further assumption on the drift, even for $p \le 1$. See again \cite{jacod2008}. \qed
\end{remark}

As noted above we have associated central limit theorems in all three cases, but for simplicity we will state the one connected to Theorem \ref{jacodlln} (a) only in the case of a continuous $X$ in which it holds irrespective of $p$. In general, such a result is expected to hold only with $p < 1$, but with additional assumptions regarding the jumps then. Similarly, the central limit theorem associated to Theorem \ref{jacodlln} (b) only holds for $p > 3$. The mode of convergence is always ($\fF$-)stable convergence in law, which means in particular that the limiting variables are typically defined on an appropriate extension of $\proba$. For details on stable convergence see Section 2.2.1 in \cite{JacPro12}. 

\begin{theorem} \label{thmjacjoint}
Let $X$ be a semimartingale of the form (\ref{ItoSemimart}).
\begin{itemize}
	\item[(a)] Suppose that $X$ is continuous and assume that Condition \ref{condclt} holds. With $f_n(x) = \Delta_n^{1-p/2} |x|^p$ we have the stable convergence 
	\[
	\De_n^{-1/2} \Big(U_n - m_p \int_0^1  \si_s^p ds\Big) \tols Y = \sqrt{m_{2p} - m_{p}^2} \int_0^1 \si_s^p dW'_s
	\]
	where $W'$ denotes an independent Brownian motion on a suitable extension of the original probability space. 
	\item[(b)] Let $p > 3$ and suppose that $X$ allows for jumps and that $X$ and $\si$ never jump at the same time. Under Condition \ref{condlln} and with $f_n(x) =  |x|^p$ for all $n$ we have the stable convergence 
	\[
	\De_n^{-1/2} \Big(U_n - \sum_{0 < s \le 1} |\De X_s|^p \Big) \tols Z = \sum_{r=1}^\infty p \text{ sign}(\De X_{S_r}) |\De X_{S_r}|^{p-1} \si_{S_r} N_r
	\]
where $(S_r)_{r \ge 1}$ denotes a sequence of stopping times exhausting the jumps of $X$ over $[0,1]$, and where $(N_r)_{r \ge 1}$ is a sequence of independent standard normal variables, also defined on a suitable extension of the original probability space. 
	\item[(c)] Suppose that $X$ allows for jumps and that $X$ and $\si$ never jump at the same time. Under Condition \ref{condlln} and with $f_n(x) =  |x|^2$ for all $n$ we have the stable convergence 
	\[
	\De_n^{-1/2} \left(U_n - [X,X]_1 \right) \tols Y + Z,
	\]
	with $Y$ as in part (a) and $Z$ as in part (b), where $W'$ and $(N_r)_{r \ge 1}$ are defined on the same extended probability space and independent.
\end{itemize}
\end{theorem}

For a proof see Theorem 5.3.6, Theorem 5.1.2 and Theorem 5.4.2 of \cite{JacPro12}.

\begin{remark}
The limiting variable in part (a) of Theorem \ref{thmjacjoint} is mixed normal with conditional variance 
\[
V = (m_{2p} - m_p^2 )\int_0^1 \si_s^{2p} ds.  
\]
Given a consistent estimator $V_n$ for $V$, Slutsky's lemma for stable convergence yields
\begin{align} \label{cltstand}
\frac{\De_n^{-1/2} \left(U_n - m_p \int_0^1 \si_s^{p} ds \right)}{\sqrt{V_n}} \tol \mathcal N(0,1).
\end{align}

In general, the central limit results connected with jumps do not allow for a mixed normal limit. An exception is the case where $\si$ and $X$ have no common jumps (compare e.g.\ Proposition 5.1.1 in \cite{JacPro12}), which is why we work under this assumption. In this case we obtain 
\[
V = \sum_{0 < s \le 1} p^2 |\De X_s|^{2p-2} \si^2_s 
\] 
for part (b) and 
\[
V = 2 \int_0^1 \si_s^{4} ds + \sum_{0 < s \le 1} 4 |\De X_s|^{2} \si^2_s .  
\]
for part (c), respectively. The goal then again is to find a consistent estimator for $V$, from which central limit theorems similar to (\ref{cltstand}) can be concluded. \qed
\end{remark} 

Historically, estimators for the asymptotic conditional variances in Theorem \ref{thmjacjoint} have been built using the exact representation of $V$ and somewhat similar statistics as the original power variations. For example, in case (a) above it is obvious from Remark \ref{remlln} that 
\[
\wV_n = \frac{m_{2p} - m_p^2}{m_{2p}} \sum_{i=1}^n g_n(\De_i^n X)  
\]
with $g_n(x) = \De_n^{1-p} |x|^{2p}$ consistently estimates $V$. In the other two cases estimation of the conditional variances is possible, yet severely more complicated due to the mixture of jumps and volatility. Plain power variations cannot be used anymore, but a truncated version where only increments $\De_i^n X$ with $|\De_i^n X| > \al \De_n^{\vpi}$, $\vpi < 1/2$, $\al > 0$, are used, combined with a local estimator for the volatility, still does the trick. See for example Theorem 9.5.1 in \cite{JacPro12}. This feature in fact is typical in high-frequency analysis: The conditional variance is often substantially more difficult to estimate than the original quantities of interest. 

\subsection{Universal estimators in the continuous case}

Two competing procedures have recently been proposed in the literature which do not try to mimic the specific structure of the limiting conditional variance, but rather construct estimators directly from the form of the original statistics $U_n$. Let us remain in the framework of power variations, so 
\[
U_n = \sum_{i=1}^n f_n(\De_i^n X),
\]
and let us write the limiting variables in Theorem \ref{jacodlln} as 
\[
U = \sum_{i=1}^n \te_{[(i-1)\De_n, i\De_n]}, 
\]
so for example 
\[
\te_{[(i-1)\De_n, i\De_n]} = [X,X]_{i \De_n} - [X,X]_{(i-1) \De_n}  = \int_{(i-1)\De_n}^{i \De_n} \si_s^2 ds + \sum_{(i-1)\De_n < s \le i\De_n} |\De X_s|^2
\]
in case of part (c). The essential idea behind the estimator from \cite{mykzha2017} is the intuition that each summand $f_n(\De_i^n X)$ within $U_n$ is in fact a local estimate for the corresponding $\te_{[(i-1)\De_n, i\De_n]}$, and this intuition remains true if several increments are aggregated. Precisely, 
\begin{align} \label{defwte}
\wte_{[i\De_n, (i+k_n)\De_n]} = \sum_{j=1}^{k_n} f_n(\De_{i+j}^n X)
\end{align}
with an auxiliary sequence $k_n \to \infty$, $k_n \De_n \to 0$, serves as an estimator for 
\begin{align} \label{defte}
\te_{[i\De_n, (i+k_n)\De_n]} = \sum_{j=1}^{k_n} \te_{[(i+j-1)\De_n, (i+j)\De_n]}
\end{align}
They therefore base their estimator on 
\begin{align*}
QV_n(k_n) = \frac{1}{k_n} \sum_{i=k_n}^{n-k_n} \left(\wte_{[(i-k_n)\De_n, i\De_n]} - \wte_{[i\De_n, (i+k_n)\De_n]}\right)^2
\end{align*}
which, using a simple decomposition, essentially mimics twice the asymptotic variance, plus an additional term due the difference of $\te_{[i\De_n, (i+k_n)\De_n]}$ and $\te_{[(i-k_n)\De_n, i\De_n]}$. When the latter approximation error is not too large compared with the other two terms, it is possible to get rid of it by working with a suitable linear combination of two different $QV_n(k_n)$. Among other possible linear combinations \cite{mykzha2017}  choose 
\[
T_n = \frac 23 \big(QV_n(k_n) - \frac 14 QV_n(2k_n) \big).
\]
An estimator for $V$ is then given by $n T_n$.  

The estimator from \cite{chrietal2017} is based on a subsampling procedure. They set 
\[
U_l^n = {k_n} \sum_{i=1}^{\lfloor \frac{n}{k_n} \rfloor} f_n(\De_{(i-1)k_n + l}^n X)
\]
for each $l = 1, \ldots, k_n$. Up to edge effects this is the same estimator as the original one, but where only each $k_n$th increment is taken into account, thus the estimator is blown up by the factor $k_n$. Again, $f_n(\De_{(i-1)k_n + l}^n X)$ is a local estimator for $\te_{[((i-1)k_n + l-1)\De_n, ((i-1)k_n + l)\De_n]}$, and if neighboring $\te_{[((i-1)k_n + l-1)\De_n, ((i-1)k_n + l)\De_n]}$ are close the each other, then $U^n_l$ should behave in the same way as the original $U_n$. In particular, a central limit theorem should hold with the same asymptotic variance, but the rate of convergence should drop to $(k_n \De_n)^{1/2}$. Therefore, the subsampling estimator for the asymtotic variance is given by
\[
\wSi_n = \frac{1}{k_n} \sum_{l=1}^{k_n} (k_n \De_n)^{-1} (U_l^n - U_n)^2,
\]
where $U_n$ serves as an approximation for the unknown limit $U$. As the convergence of $U_n$ to $U$ happens at a faster rate than the convergence of $U_l^n$ to $U$, this replacement does not cause any troubles in the limit. 

Both estimators, $n T_n$ and $\wSi_n$, are known to work in a variety of situations if $k_n \to \infty$ and $k_n \De_n \to 0$ hold and are by no means restricted to power variations. \cite{mykzha2017} work with a structural assumption and show that their estimator works in most cases where the limiting variable takes the form 
\[
U = \int_0^1 \te_s ds
\]
for some semimartingale $\te$, whereas \cite{chrietal2017} establish consistency of their subsampling estimators explicitly for power and bipower variations, including a truncated version when additional jumps are present in the process and a pre-averaged version when the process is only observed with noise. In particular, in both papers the case of a limit governed by jumps is excluded, intuitively because the implicit assumption fails that estimators close nearby will estimate the same quantity. In fact, they estimate very different quantities if a jump is present because it falls into just one interval and not into the next one.

\begin{example} \label{exfail}
Suppose that $X_t = \si W_t + J_t$ for a constant $\si > 0$ and a Poisson process $J$ with parameter $\la > 0$. Then, with 
\[
f_n(\De_i^n X) = |X_{i \De_n} - X_{(i-1)\De_n}|^2 
\]
and 
\[
U_n = \sum_{i=1}^n f_n(\De_i^n X),
\]
we have 
\[
\De_n^{-1/2} \left(U_n - [X,X]_1 \right) \tols Y + Z,
\]
according to Theorem \ref{thmjacjoint}, where the limiting variance is given by 
\[
V = 2 \si^4 +  4 \si^2 J_1.  
\] 
But, for any choice of $k_n \to \infty$ and $k_n \De_n \to 0$ we neither have $n T_n \pn V$ nor $\wSi_n \pn V$. A proof of this result will be given in the Appendix.
\end{example}

\subsection{Three new universal estimators}

In order to circumvent the problem that a jump falls into just one interval, we will present several novel estimators in the following, all of which are based on the following intuition: We fix a local interval $[i\De_n, (i+k_n)\De_n]$ first, and we will always compare two estimators constructed from increments within this interval only. These estimators are defined in such a way that a possible jump dominates both estimators in the same way, so that it is wiped out to first order. Afterwards, the local estimators based on $[i\De_n, (i+k_n)\De_n]$ are aggregated into a global estimator. 

This procedure is explained easiest for a first estimator $V_n$ which is not universal in the sense that $V_n \pn V$ holds in all three cases. Recall (\ref{defwte}) and (\ref{defte}). We will use $\wte_{[i\De_n, (i+k_n)\De_n]}$ as a  local estimator for $\te_{[i\De_n, (i+k_n)\De_n]}$ again, but it will be compared with a local power variation based on the increment $X_{(i+k_n) \De_n} - X_{i\De_n}$ which, using the same $p > 0$, also is a local estimator for $\te_{[i\De_n, (i+k_n)\De_n]}$. Recall that a possible scaling depends on the length of the interval over which the increment is computed, so the factor will be based on $k_n \De_n$ instead of $\De_n$. For example, in the continuous case we set
\[
U^n_{[i\De_n, (i+k_n)\De_n]} = (k_n \Delta_n)^{1-p/2} |X_{(i+k_n)\De_n} - X_{i\De_n}|^p
\]  
and otherwise 
\[
U^n_{[i\De_n, (i+k_n)\De_n]} = |X_{(i+k_n)\De_n} - X_{i\De_n}|^p.
\]  
The first estimator is then given by
\[
\wV_n = \frac{n}{k_n (k_n-1)} \sum_{i=0}^{n-k_n} \left(U^n_{[i\De_n, (i+k_n)\De_n]} - \wte_{[i\De_n, (i+k_n)\De_n]} \right)^2.
\]

\begin{theorem} \label{thmpowerjumps}
Let $X$ be of the form (\ref{ItoSemimart}) and let $k_n \to \infty$ such that $k_n = o(n)$. 
\begin{itemize}
	\item[(a)] Suppose that $X$ is continuous and assume that Condition \ref{condclt} holds. We have 
	\[
	\wV_n \pn V = (m_{2p} - m_p^2 )\int_0^1 \si_s^{2p} ds.  
	\]
	\item[(b)] Let $p > 3$ and suppose that $X$ allows for jumps and that $X$ and $\si$ never jump at the same time. Under Condition \ref{condlln} and with $f_n(x) =  |x|^p$ for all $n$ we have the stable convergence 
	\[
	\wV_n \tols V^* = \sum_{r=1}^\infty p^2 |\De X_{S_r}|^{2p-2} \si_{S_r}^2 (1+R_r)
	\]
	where $(S_r)_{r \ge 1}$ denotes a sequence of stopping times exhausting the jumps of $X$ over $[0,1]$ and where $(R_r)_{r \ge 1}$ denotes a sequence of i.i.d.\ random variables, independent of $\fF$ and defined on a suitable extension of the original probability space. The random variables $R_r$ have mean zero and variance one and are bounded from below by $-1$. 
	\item[(c)] Suppose that $X$ allows for jumps and that $X$ and $\si$ never jump at the same time. Under Condition \ref{condlln} we have with $f_n(x) =  |x|^2$ for all $n$
	\[
	\wV_n \tols 2 \int_0^1 \si_s^{4} ds + \sum_{r=1}^\infty 4 |\De X_{S_r}|^2 \si_{S_r}^2 (1+R_r),
	\]
	with $(S_r)_{r \ge 1}$ and $(R_r)_{r \ge 1}$  as in (b).
\end{itemize}
\end{theorem}
 
\begin{remark} \label{remvn1}
Let us discuss the heuristics behind Theorem \ref{thmpowerjumps} by distinguishing the two cases of $X$ being continuous and $X$ having jumps. The mixed case typically just combines those arguments. 
\begin{itemize}
	\item[(i)] In the continuous case, let us discuss the related, asymptotically equivalent, estimator 
\[
\wV_n^{(1)} = \frac{1}{k_n} \sum_{\ell=0}^{k_n-1} \frac{n}{k_n} \sum_{i=0}^{\lfloor \frac n{k_n} \rfloor-1} \left(U^n_{[(ik_n+\ell)\De_n, ((i+1)k_n +\ell)\De_n]} - \wte_{[(ik_n+\ell)\De_n, ((i+1)k_n +\ell)\De_n]}\right)^2
\] 
which is the same as $\wV_n$ up to small order edge effects. Note that for each fixed $\ell$ the estimator is based on observations from non-overlapping intervals. Later on these are aggregated in some type of sample mean. Then, if we set
	\[
	\te_{[u,v]} = m_p \int_u^v \si_s^p ds
	\] 
	for $u < v$, following the same proof as Theorem \ref{thmjacjoint} (a), it is easy to see that 
	\[
	\wte_{[(ik_n+\ell)\De_n, ((i+1)k_n +\ell)\De_n]} - \te_{[(ik_n+\ell)\De_n, ((i+1)k_n +\ell)\De_n]} = o_\P(k_n \De_n),
	\]
	uniformly in $i$ and $\ell$. Young's inequality allows us to replace one term by the other. As we work over disjoint intervals, 
	we then use the intuition that the
	\[
	\sqrt{\frac{n}{k_n}} \sum_{i=0}^{\lfloor \frac n{k_n} \rfloor-1} \left(U^n_{[(ik_n+\ell)\De_n, ((i+1)k_n +\ell)\De_n]} - \te_{[(ik_n+\ell)\De_n, ((i+1)k_n +\ell)\De_n]} \right)
	\]
	obey the same central limit theorem as Theorem \ref{thmjacjoint} (a). In particular, using conditional independence, it is no surprise that each 
	\[
	\frac{n}{k_n} \sum_{i=0}^{\lfloor \frac n{k_n} \rfloor-1} \left(U^n_{[(ik_n+\ell)\De_n, ((i+1)k_n +\ell)\De_n]} - \wte_{[(ik_n+\ell)\De_n, ((i+1)k_n +\ell)\De_n]}\right)^2
	\]
	estimates $V$. So does $\wV_n^{(1)}$.
	\item[(ii)] Whenever jumps are present, the idea is to implicitly assume that there are only finitely many of them and that each interval $(i\De_n, (i+k_n)\De_n]$ contains either no jump or exactly one jump. The proof of Theorem \ref{thmjacjoint} (b) shows, due to $p > 3$, that only those intervals with jumps play a role to first order in the asymptotics. For each jump time $S_r$ and for each interval such that $S_r \in (i\De_n, (i+k_n)\De_n]$, a Taylor expansion gives
	\begin{align*}
	&\wte_{[i\De_n, (i+k_n)\De_n]} - \te_{[i\De_n, (i+k_n)\De_n]} \\ =& p \text{ sign}(\De X_{S_r}) |\De X_{S_r}|^{p-1} \si_{S_r} (W_{(i+k_n)\De_n} - W_{i \De_n}) + o_\P((k_n \De_n)^{1/2})
	\end{align*}
	uniformly in $i$, where $\te_{[i\De_n, (i+k_n)\De_n]} = |\De X_{S_r}|^{p}$ and by using that $\si$ is continuous at $S_r$ by assumption. Therefore
\[
\wV_n = \sum_{r} p^2 |\De X_{S_r}|^{2p-2} \si^2_{S_r} \frac{n}{k_n^2}  \sum_{j=1}^{k_n} (W_{(i_r+k_n-j)\De_n} - W_{(i_r-j) \De_n})^2 \big(1+ o_\P(1)\big),
\]
where $((i_r-1) \De_n, i_r\De_n]$ denotes the interval which includes $S_r$. Note that the second sum above consists of highly correlated Brownian increments, and it is easy to see that its expectation and its variance are both equal to one, at least to first order. This explains the properties of the limiting distribution. \qed
\end{itemize}
\end{remark}

The lesson told by Remark \ref{remvn1} is that we need less dependence between the Brownian increments over those intervals where jumps are detected. A natural second statistic therefore is given by
\begin{align*}
\WV_n &= \frac{n}{2} \sum_{i=0}^{n-k_n} \bigg(\frac{1}{\binom{k_n}{2}}\sum_{i < u < v \le i+k_n} \big(f_n(\De_u^n X + \De_v^n X) - (f_n(\De_u^n X) + f_n(\De_v^n X))\big)^2 \bigg),
\end{align*}
where the scaling within $f_n$ again depends on the length of the corresponding interval. Let us explain the main idea behind $\WV_n$ by using the simplifying assumption again that there are only finitely many jumps which are separated in the sense that no interval $(i\De_n, (i+k_n)\De_n]$ contains more than one jump. Then in the jump case 
\[
\WV_n = \sum_{r} \frac n{k_n(k_n-1)} \sum_{j=1}^{k_n} \bigg(\sum_{\substack{i_r-j < v \le i_r+k_n-j \\ v\neq i_r}} \big(|\De_{i_r}^n X + \De_v^n X|^p - |\De_{i_r}^n X|^p \big)^2 \bigg) \big(1+ o_\P(1)\big),
\]
as only the cases with $u = i_r$ or $v = i_r$ give dominating terms to first order. If one now uses a Taylor expansion and keeps $j$ fixed first, we obtain 
\[
\WV_n = \sum_{r} \frac 1{k_n} \sum_{j=1}^{k_n} \frac n{k_n-1} \sum_{\substack{i_r-j < v \le i_r+k_n-j \\ v\neq i_r}} p^2 (\De X_{S_r})^{2p-2} \si^2_{S_r} (\De_v^n W)^2  + o_\P(1),
\]
and it is clear that we have indeed convergence in probability to the correct quantity. 

The drawback, however, is that the statistic does not converge in probability to the correct variance if the continuous part dominates. The reason is simple: We now subtract $(f_n(\De_u^n X) + f_n(\De_v^n X))$ only which is just a sum of two terms. Previously, when discussing $\wV_n$, we subtracted a sum of $k_n$ terms which asymptotically equals a functional of $\si^p$. This allowed us to mimic the arguments from the original central limit theorem. Now we estimate a quantity which is in general different from $V$. A remarkable exception is the case $p=2$ where we exactly estimate the variance $V$.

\begin{theorem} \label{thmpowerjumps2}
Let $X$ be of the form (\ref{ItoSemimart}) and let $k_n \to \infty$ such that $k_n = o(n)$. 
\begin{itemize}
	\item[(a)] Suppose that $X$ is continuous and assume that Condition \ref{condclt} holds. With 
	\[
	c_p = 2 \E \bigg[\Big( \big|\frac 1{\sqrt 2}(N_{1} + N_{2})\big|^p - \frac 12 (|N_{1}|^p + |N_{2}|^p)\Big)^2\bigg] 
	\]
	for independent standard normal $N_1$, $N_2$ we have
	\[
	\WV_n \pn c_p \int_0^1 \si_s^{2p} ds.
	\]
	\item[(b)] Let $p > 3$ and suppose that $X$ allows for jumps and that $X$ and $\si$ never jump at the same time. Under Condition \ref{condlln} and with $f_n(x) =  |x|^p$ for all $n$ we have 
	\[
	\WV_n \pn V = \sum_{0 < s \le 1} p^2 |\De X_s|^{2p-2} \si^2_s.
	\]
	\item[(c)] Suppose that $X$ allows for jumps and that $X$ and $\si$ never jump at the same time. Under Condition \ref{condlln} and with $f_n(x) =  |x|^2$ for all $n$ we have 
	\[
	\WV_n \pn V= 2 \int_0^1 \si_s^{4} ds + \sum_{0 < s \le 1} 4 |\De X_s|^{2} \si^2_s.
	\]
\end{itemize}
\end{theorem}

\begin{remark}
Note that Theorem \ref{thmpowerjumps2} (c) proves that $\WV_n$ is a consistent estimator for the asymptotic conditional variance when the quadratic variation is to be estimated. In this situation various estimators are known in the literature which all mimic the specific form of the variance; see for example Chapter 9.5 in \cite{JacPro12} or \cite{veraart2010}.  \qed
\end{remark} 

The construction of a universal estimator which converges in probability to $V$ in all three cases now combines the best from both worlds. Let $\ell_n \to \infty$ with $\ell_n = o(k_n)$ be another auxiliary sequence and set 
\begin{align*}
V_n &= \frac{n}{\ell_n(\ell_n-1)} \sum_{i=0}^{n-k_n} \frac{1}{\binom{k_n}{\ell_n}}\sum_{i \le j_1 < \ldots < j_{\ell_n} \le i+k_n} \bigg( f_n \Big(\sum_{m=1}^{\ell_n} \De_{j_m}^n X\Big) - \sum_{m=1}^{\ell_n} f_n(\De_{j_m}^n X)\bigg)^2.
\end{align*}
We see that a jump in $\De_{j_1}^n X$, say, comes together with a growing number of increments which are sufficiently independent from each other in order to ensure convergence in probability as for $\WV_n$. Also, as we subtract $\sum_{m=1}^{\ell_n} f_n(\De_{j_m}^n X)$, we consistently estimate a local version of $\si^p$ in the continuous case. Note that $\wV_n$ and $\WV_n$ are special cases with $\ell_n = k_n$ and $\ell_n =2$, respectively.  

\begin{theorem} \label{thmpowerjumps3}
Let $X$ be of the form (\ref{ItoSemimart}) and let $\ell_n, k_n \to \infty$ with $\ell_n = o(k_n)$ and $k_n = o(n)$. 
\begin{itemize}
	\item[(a)] Suppose that $X$ is continuous and assume that Condition \ref{condclt} holds. We have 
	\[
	V_n \pn V= (m_{2p} - m_p^2) \int_0^1 \si_s^{2p} ds.
	\] 
	\item[(b)] Let $p > 3$ and suppose that $X$ allows for jumps and that $X$ and $\si$ never jump at the same time. Under Condition \ref{condlln} and with $f_n(x) =  |x|^p$ for all $n$ we have 
	\[
	V_n \pn V = \sum_{0 < s \le 1} p^2 |\De X_s|^{2p-2} \si^2_s.
	\]
	\item[(c)] Suppose that $X$ allows for jumps and that $X$ and $\si$ never jump at the same time. Under Condition \ref{condlln} and with $f_n(x) =  |x|^2$ for all $n$ we have 
	\[
	V_n \pn V= 2 \int_0^1 \si_s^{4} ds + \sum_{0 < s \le 1} 4 |\De X_s|^{2} \si^2_s.
	\]
\end{itemize}
\end{theorem}

\section{Conclusion} 
\def\theequation{4.\arabic{equation}}
\setcounter{equation}{0}

In this paper we have presented a new class of estimators for the asymptotic (conditional) variance in limit theorems for semimartingales. These estimators are only based on the form of the original statistics 
\[
U_n = \sum_{i=1}^n f_n(\De_i^n X)
\]
in the central limit theorem, and we have shown in Theorem \ref{thmpowerjumps3} that they are consistent for power variations in all three possible regimes: For a dominating continuous martingale part, for dominating jumps and for the quadratic variation. 

Even though the estimator $V_n$ discussed in Theorem \ref{thmpowerjumps3} gives a positive answer to the question whether such universal estimators exist, its application in practice is difficult, as we need to compute statistics over each of the $\binom{k_n}{\ell_n}$ subintervals within $(i\De_n, (i+k_n)\De_n]$ in order to obtain $V_n$. From a computational point of view this is certainly not a reasonable strategy, at least under the conditions $\ell_n \to \infty$ and $\ell_n = o(k_n)$. The other estimators $\wV_n$ and $\WV_n$ are constructed with $\ell_n = k_n$ and $\ell_n = 2$, respectively, so they are computationally much less expensive, though not consistent in all situations. 

Future research clearly needs to investigate the practical properties of this new class of estimators, for $\wV_n$ in comparison to \cite{mykzha2017} and \cite{chrietal2017} in the continuous case, but also with a focus towards the properties of $\WV_n$ in the case of quadratic variation. This new estimator is consistent in all situations, with jumps or not, so one does not need to test in advance whether jumps are present in the path of $X$ or not.

\section{Proofs} \label{sec:proof}
\def\theequation{5.\arabic{equation}}
\setcounter{equation}{0}

Throughout the proofs we will assume that the processes $(b_s)$, $(\sigma_s)$ and $(X_s)$ are bounded, and we will also assume that $|\delta(s,z)|$ is bounded by a deterministic function $\gamma(z)$ satisfying $\int \ga^2(z) \la(dz) < \infty$. In fact, according to Condition \ref{condlln} we know that $(b_s)$ and $(\delta(s,z))$ safisfy such claims locally, and we also know that $(\si_s)$ is c\`adla\`g, and then a standard localization procedure as in Section 4.4.1 in \cite{JacPro12} shows that we may assume global bounds without loss of generality. Similarly, whenever we explicitly need Condition \ref{condclt}, we may further assume that $(\Wsi_s)$, $(\Wb_s)$, $(\al_s)$ and $(\al'_s)$ are bounded as well, and we may also assume that $(\si)$ is bounded away from zero. Also, $C > 0$ denotes a universal constant which may change from line to line, and we write $C_r$ whenever we want to emphasize dependence of the constant on an auxiliary parameter such as $r$.

We introduce the decomposition $X_t=X_0+B(q)_t+X^c_t+M(q)_t+N(q)_t$ of the It\^o semimartingale \eqref{ItoSemimart} with
\begin{align*}
B(q)_t&=\int_0^t \Big( b_s -\int(\delta(s,z)\mathds{1}_{\{|\delta(s,z)|\leq 1\}}-\delta(s,z)\mathds{1}_{\{\gamma(z)\leq1/q\}})\lambda(dz)\Big)ds,
\\X^c_t &= \int_0^t \sigma_s dW_s,
\\M(q)_t&=\int_0^t \int \delta(s,z) \mathds{1}_{\{\gamma(z)\leq 1/q\}}(\mu-\nu)(ds,dz),
\\N(q)_t&=\int_0^t \int \delta(s,z) \mathds{1}_{\{\gamma(z)>1/q\}}\mu(ds,dz).
\end{align*}
Here $q$ is a parameter which controls whether jumps are classified as small jumps or big jumps. We also set $X(q)_t = B(q)_t + X^c_t + M(q)_t$ and denote the derivative process of $B(q)$ with $b(q)$. From the integrability condition on $\ga$ one immediately obtains $|b(q)| \le Cq$.

\subsection{Proof of Example \ref{exfail}}
Let $A$ be the subset of $\Om$ such that $J$ contains exactly one jump in $(0,1)$ and that the jump time $S$ is in $(0,1) \backslash \Q$. Obviously, $\P(A) > 0$, and it is sufficient to prove that both $n T_n \mathds{1}_A$ and $\wSi \mathds{1}_A$ diverge to infinity in probability.

For $T_n$, on $A$, suppose that $n$ is large enough such that $k_n \De_n < S < 1-k_n \De_n$. Then, each
\[
QV_n(k_n) = \frac{1}{k_n} \sum_{i=k_n}^{n-k_n} (\wte_{[(i-k_n)\De_n, i\De_n]} - \wte_{[i\De_n, (i+k_n)\De_n]})^2
\]
consists of $2k_n$ summands which are affected by the one jump and of $n-4k_n+1$ summands which are not. Suppose for example that $i = \left\lceil nS  \right\rceil$. Then 
\begin{align*}
& \wte_{[(i-k_n)\De_n, i\De_n]} - \wte_{[i\De_n, (i+k_n)\De_n]} \\
=& \sum_{j=1}^{k_n-1} \si^2 (|\De_{i-k_n+j}^n W|^2 - |\De_{i+j}^n W|^2) + |1+\si \De_i^n W|^2 - \si^2  |\De_{i+k_n}^n W|^2 \\
=& 1+ 2 \si \De_i^n W + \sum_{j=1}^{k_n} \si^2 (|\De_{i-k_n+j}^n W|^2 - |\De_{i+j}^n W|^2) = 1 + O_\P(\sqrt{\De_n}), 
\end{align*}
where we have used $k_n \De_n \to 0$. Consequently, 
\begin{align*}
\frac{1}{k_n} \sum_{i=\left\lceil S \De_n^{-1} \right\rceil- k_n}^{\left\lceil S \De_n^{-1} \right\rceil + k_n -1} (\wte_{[(i-k_n)\De_n, i\De_n]} - \wte_{[i\De_n, (i+k_n)\De_n]})^2 = 2 + O_\P(\sqrt{\De_n}).
\end{align*}
The sum over the remaining $n-4k_n+1$ terms asymptotically behaves in the same way as the entire $QV_n(k_n)$ in the case without jumps and is of order $\De_n$ according to Theorem 4 of \cite{mykzha2017}. Therefore 
\[
T_n = \frac 23 \big(QV_n(k_n) - \frac 14 QV_n(2k_n) \big) = 1 + O_\P(\sqrt{\De_n}),
\] 
and $n T_n$ diverges on $A$. 

Similarly, on the set $A$ we have that only one of the statistics $U_l^n$ contains the increment with the one jump, whereas the remaining $k_n-1$ intervals are not affected by it. Therefore, each of the latter statistics satisfies $U_l^n - U_n = O_\P(1)$ as restricted to $A$ both statistics converge in probability to $\si^2$ and $\si^2 + 1$, respectively. We conclude that 
\[
\wSi_n = \frac{1}{k_n} \sum_{l=1}^{k_n} (k_n \De_n)^{-1} (U_l^n - U_n)^2 = O_\P((k_n \De_n)^{-1})
\]
on $A$, so it does not converge as well. \qed

\subsection{Proof of Theorems \ref{thmpowerjumps}, \ref{thmpowerjumps2} and \ref{thmpowerjumps3}} 
We will proceed as follows: In all cases we will only show parts (a) and (b), and we will discuss these in separate sections. The proof of part (c) mostly just combines the ideas from (a) and (b) after one separates intervals with and without jumps of $N(q)$. Within each section we will start with the result from Theorem \ref{thmpowerjumps3} which we will prove in essentially all details. Afterwards we discuss the necessary changes for Theorems \ref{thmpowerjumps} and \ref{thmpowerjumps2}. Note that we can use analogous proofs for most parts because the estimators are essentially all the same, just with $\ell_n$ varying between 2 and $k_n$.

Before we begin with the proofs of the main theorems, we provide a key lemma which will be used extremely often throughout the remaining sections. 
\begin{lemma} \label{lemyoung}
Let 
\[ 
X_n = \sum_{i = 1}^{n-k_n} (\chi_i^n)^2
\]
and suppose that there exists 
\[
R_n = \sum_{i = 1}^{n-k_n} (\rho_i^n)^2
\]
such that $R_n \weak X$ and 
\begin{align} \label{rmx}
\sum_{i = 1}^{n-k_n} (\chi_i^n - \rho_i^n)^2 \pn 0.
\end{align}
Then $X_n \weak X$.
\end{lemma}

\textbf{Proof:} We will only show $X_n - R_n \pn 0$. Note that for each $\eps>0$ there exists some $C_\eps > 0$ such that 
\begin{align} \label{ineqstd}
|(x+y)^2 - x^2| \le \eps x^2 + C_\eps y^2,
\end{align}
which is a simple consequence of Young's inequality. Therefore 
\begin{align*}
|X_n - R_n| \le \eps R_n + C_\eps \sum_{i = 1}^{n-k_n} (\chi_i^n - \rho_i^n)^2,
\end{align*}
and we obtain
\begin{align*}
\limsup_{n \to \infty} \P(|X_n - R_n| > \de) \le \limsup_{n \to \infty} \P \Big(R_n \ge \frac{\de}{2\eps} \Big) \le \P \Big(X \ge \frac{\de}{2\eps} \Big)
\end{align*}
for each fixed $\eps$, where we have first used (\ref{rmx}) and the Portmanteau theorem plus $R_n \weak X$ afterwards. Letting $\eps \to 0$ then finishes the proof. \qed

\subsubsection{Proof of part (a)}
We will start with Theorem \ref{thmpowerjumps3} and discuss $V_n$. In the situation of a continuous $X$ a simple computation using the respective standardisation of $f_n$ shows that the estimator reads as 
$$
V_n = \frac {\ell_n}{n (\ell_n-1)} \sum_{i=0}^{n-k_n} \frac{1}{\binom{k_n}{\ell_n}} \sum_{1 \le j_1 < \ldots < j_{\ell_n} \le k_n} \left(V_{i+j_1, \ldots, i+j_{\ell_n}}^n \right)^2
$$
with 
\[
V_{i+j_1, \ldots, i+j_{\ell_n}}^n =  (\ell_n \De_n)^{-p/2} \Big|\sum_{m=1}^{\ell_n} \De_{i+j_m}^n X \Big|^p - \frac 1{\ell_n} \sum_{m = 1}^{\ell_n} \De_n^{-p/2} |\De_{i + j_m}^n X|^p.
\]

The main strategy in the proof of $V_n \pn V$ is to apply Lemma \ref{lemyoung} several times, which means that one successively replaces $V_n$ by simpler terms until one ends up with 
\[
\overline V_n = \frac{1}{n} \sum_{i=0}^{n-k_n} \si^{2p}_{i \De_n} \frac{1}{\binom{k_n}{\ell_n}} \sum_{1 \le j_1 < \ldots < j_{\ell_n} \le k_n} \bigg( (\ell_n \De_n)^{-p/2} \Big|\sum_{m=1}^{\ell_n} \De_{i+j_m}^n W\Big|^p - m_p \bigg)^2.
\]
We first prove $\overline V_n \pn V$ for which we set $h_n(x_1, \ldots, x_{\ell_n}) = |\ell_n^{-1/2} (x_1 + \ldots + x_{\ell_n})|^p$ and 
\[
U_i^n = \frac{1}{\binom{k_n}{\ell_n}} \sum_{1 \le j_1 < \ldots < j_{\ell_n} \le k_n} \left( h_n(\De_n^{-1/2} \De_{i+j_1}^n W, \ldots, \De_n^{-1/2} \De_{i+j_{\ell_n}}^n W) - m_p \right)^2.
\]
Clearly, $\E[U_i^n] = m_{2p} - m_p^2$, and because of conditional independence, boundedness of $(\si)$ and the Cauchy-Schwarz inequality we also have
\begin{align} \label{varU}
& \E \bigg[\Big(\frac{1}{n} \sum_{i=0}^{n-k_n} \si^{2p}_{i \De_n} (U_i^n - \E[U_i^n]) \Big)^2 \bigg]  \le \frac{C}{n^2} \sum_{i,r=0}^{n-k_n} \Ind_{\{|i-r| \le k_n\}} \sqrt{ \Var(U_i^n) \Var(U_r^n) }. 
\end{align}
Using Theorem 1.2.3 in \cite{denker1985} on an upper bound for the variance of a U statistic we obtain
\[
\Var(U_i^n) \le \frac{\ell_n}{k_n} \Var\left((h_n(\De_n^{-1/2} \De_{i+j_1}^n W, \ldots, \De_n^{-1/2} \De_{i+j_{\ell_n}}^n W) - m_p)^2\right) \le C \frac{\ell_n}{k_n},
\]
so as a consequence of $\ell_n \De_n \to 0$ 
\[
\overline V_n = \frac{1}{n} \sum_{i=0}^{n-k_n} \si^{2p}_{i \De_n} U_i^n = \frac{1}{n} \sum_{i=0}^{n-k_n} \si^2_{i \De_n} (m_{2p} - m_p^2) + o_\P(1).
\]
Convergence in probability of the latter quantity to $V$ is standard. 

It remains to prove that the simplification to $\overline V_n$ is adequate. We first show $\underline V_n \pn V$ for
\[
\underline V_n = \frac {\ell_n}{n (\ell_n-1)} \sum_{i=0}^{n-k_n} \si^{2p}_{i \De_n} \frac{1}{\binom{k_n}{\ell_n}} \sum_{1 \le j_1 < \ldots < j_{\ell_n} \le k_n} \left(\underline V_{i+j_1, \ldots, i+j_{\ell_n}}^n\right)^2
\]
with 
\[
\underline V_{i+j_1, \ldots, i+j_{\ell_n}}^n = (\ell_n \De_n)^{-p/2} \Big|\sum_{m=1}^{\ell_n} \De_{i+j_m}^n W \Big|^p - \frac{1}{\ell_n} \sum_{m=1}^{\ell_n} \De_n^{-p/2} |\De_{i+j_m}^n W|^{p}.
\]
Using Lemma \ref{lemyoung}, boundedness of $(\si_s)$ and $\overline V_n \pn V$ we just have to establish
\begin{align} \label{hilf1}
\frac{1}{n} \sum_{i=0}^{n-k_n} \frac{1}{\binom{k_n}{\ell_n}} \sum_{1 \le j_1 < \ldots < j_{\ell_n} \le k_n} \left(\frac{1}{\ell_n} \sum_{m=1}^{\ell_n} \De_n^{-p/2} |\De_{i+j_m}^n W|^{p} - m_p\right)^2 \pn 0
\end{align}
in order to show $\frac{\ell_n-1}{\ell_n}\underline V_n \pn V$, and the claim regarding $\underline V_n$ then follows from $\ell_n \to \infty$. Note that 
\[
\overline T_{i+j_1, \ldots, i+j_{\ell_n}}^n = \left(\frac{1}{\ell_n} \sum_{m=1}^{\ell_n} \De_n^{-p/2} |\De_{i+j_m}^n W|^{p} - m_p\right)^2
\]
satisfies $\E[|\overline T_{i+j_1, \ldots, i+j_{\ell_n}}^n|^2] \le C/\ell_n$ by independence of the Brownian increments, 
so that (\ref{hilf1}) follows from 
\begin{align*}
\E \left[ \frac{1}{n} \sum_{i=0}^{n-k_n} \frac{1}{\binom{k_n}{\ell_n}} \sum_{1 \le j_1 < \ldots < j_{\ell_n} \le k_n} \left(\frac{1}{\ell_n} \sum_{m=1}^{\ell_n} \De_n^{-p/2} |\De_{i+j_m}^n W|^{p} - m_p\right)^2 \right] \le \frac{C}{\ell_n} \to 0.
\end{align*}

Finally, another application of Lemma \ref{lemyoung} together with $\underline V_n \pn V$, plus the obvious $(v+w)^2 \le 2(v^2 + w^2)$, shows that the proof of $V_n \pn V$ boils down to showing 
\begin{align} \label{contsecond}
\frac {\ell_n}{n (\ell_n-1)} \sum_{i=0}^{n-k_n} \frac{1}{\binom{k_n}{\ell_n}} \sum_{1 \le j_1 < \ldots < j_{\ell_n} \le k_n} \De_n^{-p} \left(\frac{1}{\ell_n} \sum_{m=1}^{\ell_n} \left( |\De_{i+j_m}^n X|^{p} - \si_{i \De_n}^p |\De_{i+j_m}^n W|^{p} \right) \right)^2 \pn 0
\end{align}
as well as 
\begin{align} 
\frac {\ell_n}{n (\ell_n-1)} \sum_{i=0}^{n-k_n} \frac{1}{\binom{k_n}{\ell_n}} \sum_{1 \le j_1 < \ldots < j_{\ell_n} \le k_n} \underline T_{i+j_1, \ldots, i+j_{\ell_n}}^n \pn 0  \label{contfirst}
\end{align}
with 
\[
\underline T_{i+j_1, \ldots, i+j_{\ell_n}}^n = (\ell_n \De_n)^{-p} \left(\big|\sum_{m=1}^{\ell_n} \De_{i+j_m}^n X\big|^p - \si_{i \De_n}^{p} \big|\sum_{m=1}^{\ell_n} \De_{i+j_m}^n W \big|^p \right)^2.
\]
The proof is similar for both claims, and we will only prove (\ref{contfirst}) in detail. 

To this end, let $\psi: \R \to \R$ be a smooth function such that 
\[
\mathds{1}_{[1,\infty)}(x) \le \psi(x) \le \mathds{1}_{[1/2,\infty)}(x),
\]
and for any $A > 0$ and $p > 0$ we set 
\[
\psi_A(x) = \psi\big(\frac{|x|}A\big), ~~ \psi'_A(x) = 1-\psi_A(x), ~~\psi_{A,p}(x) = \psi_A(x) |x|^p, ~~ \psi'_{A,p}(x) = \psi'_A(x) |x|^p.
\]
Clearly, 
\begin{align*}
\underline T_{i+j_1, \ldots, i+j_{\ell_n}}^n \le 2 (\underline T_{i+j_1, \ldots, i+j_{\ell_n},A}^n + \underline T_{i+j_1, \ldots, i+j_{\ell_n},A}^{'n})
\end{align*}
with 
\begin{align*}
\underline T_{i+j_1, \ldots, i+j_{\ell_n},A}^n = \left(\psi_{A,p} \left(\frac{\sum_{m=1}^{\ell_n} \De_{i+j_m}^n X}{\sqrt{\ell_n \De_n}} \right) - \psi_{A,p} \left(\si_{i \De_n} \frac{\sum_{m=1}^{\ell_n} \De_{i+j_m}^n W}{\sqrt{\ell_n \De_n}} \right)\right)^2 
\end{align*}
and similarly for $\underline T_{i+j_1, \ldots, i+j_{\ell_n},A}^{'n}$, but with $\psi_{A,p}$ replaced by $\psi'_{A,p}$. (\ref{contfirst}) then follows from
\[
\lim_{A \to \infty} \limsup_{n \to \infty} \frac {\ell_n}{n (\ell_n-1)} \sum_{i=0}^{n-k_n} \frac{1}{\binom{k_n}{\ell_n}} \sum_{1 \le j_1 < \ldots < j_{\ell_n} \le k_n}\E \left[ \underline T_{i+j_1, \ldots, i+j_{\ell_n},A}^n\right] = 0
\]
and 
\[
\lim_{n \to \infty} \frac {\ell_n}{n (\ell_n-1)} \sum_{i=0}^{n-k_n} \frac{1}{\binom{k_n}{\ell_n}} \sum_{1 \le j_1 < \ldots < j_{\ell_n} \le k_n}\E \left[ \underline T_{i+j_1, \ldots, i+j_{\ell_n},A}^{'n}\right] = 0
\]
for every fixed $A > 0$. The first claim can be quickly deduced from  
\[
\psi_{A,p}(x) = \psi \big(\frac xA \big) |x|^p \le \Ind_{\{2|x| \ge A \}} |x|^p \le \frac{2|x|^{p+1}}{A},
\]
$(v+w)^2 \le 2(v^2 + w^2)$ and e.g.\
\[
\E \left[\left| \frac{\sum_{m=1}^{\ell_n} \De_{i+j_m}^n X}{\sqrt{\ell_n \De_n}} \right|^{2p+2} \right] \le C
\]
which is a consequence of the Burkholder-Davis-Gundy inequality and the boundedness assumption for $(b_s)$ and $(\si_s)$.

So let finally be $A$ fixed. It is easy to see that $\psi'_{A,p}$ is bounded and uniformly continuous,
and it follows that 
\[
\te(\eps) = \sup_{x \in \R, |y| \le \eps} \left| \psi'_{A,p}(x+y) - \psi'_{A,p}(x) \right| \to 0
\] 
as $\eps \to 0$. In particular, 
\[
\left| \psi'_{A,p}(x+y) - \psi'_{A,p}(x) \right| \le \te(\eps) + \left| \psi'_{A,p}(x+y) - \psi'_{A,p}(x) \right| \Ind_{\{|y| > \eps\}} \le \te(\eps) + C_A \frac{y^2}{\eps^2}.
\]
By first letting $n \to \infty$ and then $\eps \to 0$ it is thus sufficient to prove 
\begin{align*} 
\frac {\ell_n}{n (\ell_n-1)} \sum_{i=0}^{n-k_n} \frac{1}{\binom{k_n}{\ell_n}} \sum_{1 \le j_1 < \ldots < j_{\ell_n} \le k_n}\E \Bigg[ \bigg(\frac{\sum_{m=1}^{\ell_n} \De_{i+j_m}^n X - \si_{i \De_n}\sum_{m=1}^{\ell_n} \De_{i+j_m}^n W}{\sqrt{\ell_n \De_n}} \bigg)^2\Bigg] \to 0
\end{align*}
as $n \to \infty$. Using $(v+w)^2 \le 2(v^2 + w^2)$ once more, we can discuss the absolutely continuous part of the increments and the Brownian parts separately, and the proof for the first terms follows from
\[
\frac 1{{\ell_n \De_n}} \E \left[\left(\int_{(i+j_1-1)\De_n}^{(i+j_1)\De_n} b_s ds + \ldots + \int_{(i+j_{\ell_n}-1)\De_n}^{(i+j_{\ell_n})\De_n} b_s ds\right)^2 \right]
\le C \ell_n \De_n \to 0.
\] 
We can thus assume $dX_t = \si_t dW_t$, and we will first prove the result in the case of a continuous $\si$. We have 
\begin{align*}
&\frac 1{{\ell_n \De_n}} \E \left[\left(\int_{(i+j_1-1)\De_n}^{(i+j_1)\De_n} (\si_s - \si_{i \De_n}) dW_s + \ldots + \int_{(i+j_{\ell_n}-1)\De_n}^{(i+j_{\ell_n})\De_n} (\si_s - \si_{i \De_n}) dW_s\right)^2\right] \\ =& \frac 1{{\ell_n \De_n}}\E \left[\left(\int_{(i+j_1-1)\De_n}^{(i+j_1)\De_n} (\si_s - \si_{i \De_n})^2 ds + \ldots + \int_{(i+j_{\ell_n}-1)\De_n}^{(i+j_{\ell_n})\De_n} (\si_s - \si_{i \De_n})^2 ds\right)^2\right],
\end{align*}
so that 
\begin{align*}
& \frac {1}{\ell_n-1} \sum_{i=0}^{n-k_n} \frac{1}{\binom{k_n}{\ell_n}} \sum_{1 \le j_1 < \ldots < j_{\ell_n} \le k_n}\E \left[\left(\sum_{m=1}^{\ell_n}\int_{(i+j_m-1)\De_n}^{(i+j_m)\De_n} (\si_s - \si_{i\De_n}) dW_s\right)^2 \right] \\ =& \frac {\ell_n}{\ell_n-1}\frac{1}{k_n} \sum_{i=0}^{n-k_n} \int_{i \De_n}^{(i+k_n) \De_n} (\si_s - \si_{i\De_n})^2 ds \le \int_{0}^{1} \frac{2}{k_n} \sum_{m=0}^{k_n-1} \E[(\si_s - \si_{(\lfloor ns \rfloor - m)^+\De_n})^2] ds
\end{align*}
where we have used that every interval $[(i+j_m-1)\De_n, (i+j_m)\De_n]$ appears $\binom{k_n-1}{\ell_n - 1}$ times and
\begin{align} \label{binkoef2}
\ell_n \binom{k_n}{\ell_n} = k_n \binom{k_n-1}{\ell_n-1},
\end{align}
plus $\ell_n \le 2(\ell_n -1)$ for any $\ell_n \ge 2$. Convergence to zero in probability then follows from continuity of $\si$ and dominated convergence.

In the general case we use the reasoning from Lemma 3.4.8 in \cite{JacPro12}. A standard argument using $\int_0^1 \sigma_s^2 ds \le C$ proves the existence of a sequence $\si(u)$ of adapted continuous processes $\si(u)$ such that 
\begin{align} \label{sigmau}
\E\left[\int_{0}^1\left( \si_s - \si(u) \right)^2 ds\right] \to 0
\end{align}
as $u \to \infty$. Thus, setting $X(u)_t = X_0 + \int_0^t b_s ds + \int_0^t \si(u)_s dW_s$, we have already shown 
\[
V_n(u) \pn V(u) = (m_{2p} - m_p^2) \int_0^1 \si(u)_s^p ds
\]
as $n \to \infty$, where $V_n(u)$ denotes the statistic $V_n$, but based on $X(u)$. Clearly, $V(u) \pn V$ as $u \to \infty$ as well, so it remains to prove 
\[
\lim_{u \to \infty} \limsup_{n \to \infty} \P(|V_n - V_n(u)| > \eta) = 0
\]
for every $\eta > 0$. Using (\ref{ineqstd}) one has to deal with similar claims as (\ref{contsecond}) and (\ref{contfirst}), but with $\si = 1$ and where $W$ becomes $X(u)$. Reproducing these lines the proof finally follows from (\ref{sigmau}).

For Theorem \ref{thmpowerjumps} the proof holds without any changes, because we have only used $\ell_n \De_n \to 0$ and $\ell_n \to \infty$ which holds for $\ell_n = k_n$ as well. The situation is different for Theorem \ref{thmpowerjumps2} in which case
\[
\underline V_n = \frac{1}{n} \sum_{i=0}^{n-k_n} \si^{2p}_{i \De_n} \underline U_i^n
\]
for a U statistic of the form
\[
\underline U_i^n = \frac{1}{\binom{k_n}{2}} \sum_{1 \le j_1 < j_{2} \le k_n} 2 \left(\underline V_{i+j_1, i+j_{2}}^n\right)^2 
\]
where
\begin{align*}
\underline V_{i+j_1, i+j_{2}}^n &= (2 \De_n)^{-p/2} |\De_{i+j_1}^n W + \De_{i+j_{2}}^n W|^p - \frac 12 \De_n^{-p/2} (|\De_{i+j_1}^n W|^{p} + |\De_{i+j_2}^n W|^{p}) \\ &= \Big|\frac 1{\sqrt 2}(N_{i+j_1} + N_{i+j_2})\Big|^p - \frac 12 (|N_{i+j_1}|^p + |N_{i+j_2}|^p)
\end{align*}
and the latter equality is to be understood in distribution, with the $N_{i+j}$ all independent standard normal. Setting 
\[
c_p = 2 \E \bigg[\bigg(\Big|\frac 1{\sqrt 2}(N_{i+j_1} + N_{i+j_2})\Big|^p - \frac 12 (|N_{i+j_1}|^p + |N_{i+j_2}|^p)\bigg)^2\bigg] 
\] 
the same reasoning as for (\ref{varU}) gives 
\[
\underline V_n = c_p \frac{1}{n} \sum_{i=0}^{n-k_n} \si^{2p}_{i \De_n} + o_\P(1) = c_p \int_0^1 \si_s^{2p} ds + o_\P(1).
\]
The remainder of the proof remains unchanged. Note finally that 
\[
	c_2 = \frac 12 \E \Big[\big(|N_1 + N_2|^2 - (|N_1|^2 + |N_{2}|^2)\big)^2\Big] = 2 \E[N_1^2 N_2^2] = 2.
\]
\qed

\subsubsection{Proof of part (b)}
We define  $L_m = \{ z ~|~ \ga(z) > 1/m\}$ for any $m \ge 1$, and let $\{S(m,j) ~|~ j \ge 1\}$ denote the jump times of the Poisson process $\mathds{1}_{L_m \backslash L_{m-1}} \star \mu$ over $[0,1]$. Then, if $(S_r)_{r \ge 1}$ is a reordering of the double sequence $(S(m,j))_{m,j \ge 1}$, we denote with $P_q$ the set of all indices $r$ such that $S_r = S(m,j)$ for some $m \le q$. By definition, these are the jump times of $N(q)$ over $[0,1]$. Further, let $\Om(n,q)$ be the set of all $\om$ on which $N(q)$ has at most one jump in each interval $[i \De_n, (i+k_n)\De_n]$, $i=0, \ldots, n-k_n$, all jumps of $N(q)$ over $[0,1]$ occur within $[k_n \De_n, 1- k_n \De_n]$ and where 
\[
|X(q)(\om)_{t+s} - X(q)(\om)_t| \le 2/q \text{ for all } t \in [0,1] \text{ and } s \in [0, k_n \De_n]. 
\] 
Since $X(q)$ is c\`adl\`ag and $N(q)$ only possesses finitely many jumps on $[0,1]$, it is clear that $\P(\Om(n,q)) \to 1$ as $n \to \infty$ for any $q > 0$. As we will typically let first $n \to \infty$ and then $q \to \infty$, we will sometimes assume $\om \in \Om(n,q)$. 

We introduce the notation $i_r$ to denote the interval $((i_r-1) \De_n, i_r \De_n]$ containing the $r$th jump $\De X_{S_r}$ of $N(q)$. In this case we have
\[
V_n = \frac {n}{\ell_n (\ell_n-1)} \sum_{i=0}^{n-k_n} \frac{1}{\binom{k_n}{\ell_n}} \sum_{1 \le j_1 < \ldots < j_{\ell_n} \le k_n} \left(\Big|\sum_{m=1}^{\ell_n} \De_{i+j_m}^n X\Big|^p - \sum_{m = 1}^{\ell_n} |\De_{i + j_m}^n X|^p \right)^2
\]
and the key to the proof will be the decomposition $V_n = V_n(q) + V_n'(q)$ with 
\[
V_n(q) = \frac {n}{\ell_n(\ell_n-1)} \sum_{r \in P_q} \sum_{\al=1}^{k_n} Y_{r,\al}^{(n)}, \quad V_n'(q) = V_n - V_n(q), 
\]
and where
\[
Y_{r,\al}^{(n)} = \frac{1}{\binom{k_n}{\ell_n}} \sum_{\substack{1 \le j_1 < \ldots < j_{\ell_n} \le k_n \\ \{j_1, \ldots, j_{\ell_n}\} \cap \{ \al \} \neq \emptyset }} \left(\Big|\sum_{m=1}^{\ell_n} \De_{i_r-\al+j_m}^n X\Big|^p - |\De_{i_r}^n X|^p \right)^2.
\]
Clearly the proof is finished once we have shown
\begin{align} \label{convnqfest}
V_n(q) \pn V(q) = \sum_{r \in P_q} p^2 |\De X_{S_r}|^{2p-2} \si^2_{S_r}
\end{align}
as $n \to \infty$ for any fixed $q$, 
\begin{align} \label{convq}
 V(q) \pn V= \sum_{0 < s \le 1} p^2 |\De X_{s}|^{2p-2} \si^2_{s}
\end{align}
as $q \to \infty$, as well as 
\begin{align} \label{stabneg}
\lim_{q \to \infty} \limsup_{n \to \infty} \P \left( \left| V_n'(q) \right|  > \eta \right) = 0
\end{align}
for all $\eta > 0$. Note that (\ref{convq}) is a direct consequence of monotone convergence. Regarding (\ref{stabneg}) we observe that increments of $X$ and $X(q)$ coincide when no jump of $N(q)$ is present. Therefore, and using $\ell_n \le 2 (\ell_n- 1)$ for $\ell_n \ge 2$, we have the inequality 
\[
|V_n'(q)| \le A_n(q) + B_n(q)
\]
with 
\[
A_n(q) = \frac {2n}{\ell_n^2} \sum_{i=0}^{n-k_n} \frac{1}{\binom{k_n}{\ell_n}} \sum_{1 \le j_1 < \ldots < j_{\ell_n} \le k_n} \left(\Big|\sum_{m=1}^{\ell_n} \De_{i+j_m}^n X(q)\Big|^p - \sum_{m = 1}^{\ell_n} |\De_{i + j_m}^n X(q)|^p \right)^2
\]
and 
\[
B_n(q) = \frac {n}{\ell_n(\ell_n-1)} \sum_{r \in P_q} \sum_{\al=1}^{k_n} |Z_{r,\al}^{(n)} - Y_{r,\al}^{(n)}|
\]
for
\[
Z_{r,\al}^{(n)} = \frac{1}{\binom{k_n}{\ell_n}} \sum_{\substack{1 \le j_1 < \ldots < j_{\ell_n} \le k_n \\ \{j_1, \ldots, j_{\ell_n}\} \cap \{ \al \} \neq \emptyset }} \left(\Big|\sum_{m=1}^{\ell_n} \De_{i_r-\al+j_m}^n X\Big|^p  -\sum_{m = 1}^{\ell_n} |\De_{i_r - \al + j_m}^n X|^p \right)^2.
\] 

We will start with the first part of \eqref{stabneg} and prove
\begin{align} \label{410}
\lim_{q \to \infty} \limsup_{n \to \infty} \P \left( A_n(q) > \eta \right) = 0,
\end{align} 
for which we set 
\begin{align} \label{defyq}
Y(q)_t = \int_{i \De_n}^t \Ind_{B_{i,j_1, \ldots, j_{\ell_n}}^n}(s) dX(q)_s, \quad t \ge i \De_n,
\end{align}
where we use the shorthand notation 
\[
B = B_{i,j_1, \ldots, j_{\ell_n}}^n = ((i+j_1-1)\De_n, (i+j_1)\De_n] \cup \ldots \cup ((i+j_{\ell_n}-1)\De_n, (i+j_{\ell_n})\De_n].
\]
We will basically apply (5.1.19) in \cite{JacPro12} which is stated for increments of $X(q)$ rather than for $Y(q)$, but the proof works similarly in our situation. Let us introduce some notation. We set $f(x) = |x|^p$ as well as 
\[
k(x,y) = f(x+y) - f(x) - f(y), \quad g(x,y) = k(x,y) - f'(x) y.
\]
Then we obtain
\begin{align*}
&|\De_{i+j_1}^n X(q) + \ldots + \De_{i+j_{\ell_n}}^n X(q)|^p = f(Y(q)_{(i+k_n)\De_n}) \\ =& \sum_{i\De_n < s \le (i+k_n)\De_n} f(\De X(q)_s) \Ind_B(s) + A(n,q,B)_{(i+k_n)\De_n} + M(n,q,B)_{(i+k_n)\De_n},
\end{align*}
where $M(n,q,B)$ is a square-integrable martingale with predictable bracket $A'(n,q,B)$, and where 
\[
A(n,q,B)_t = \int_{i\De_n}^{t} a(n,q,B)_u du, \quad A'(n,q,B) = \int_{i\De_n}^{t} a'(n,q,B)_u du,
\]
with 
\begin{align*}
a(n,q,B)_u &= f'(Y(q)_{u-}) b(q)_u \Ind_B(u) + \frac 12 f''(Y(q)_{u-}) \si^2_u \Ind_B(u)\\ &+ \int g(Y(q)_{u-}, \de(u,z)) \Ind_{\{\ga(z) \le 1/q\}} \Ind_B(u) \la(dz)
\end{align*}
and
\[
a'(n,q,B)_u = (f'(Y(q)_{u-}))^2 \si^2_u \Ind_B(u) +  \int k(Y(q)_{u-}, \de(u,z))^2 \Ind_{\{\ga(z) \le 1/q\}} \Ind_B(u) \la(dz). 
\]
Similarly, 
\begin{align*}
&\sum_{m = 1}^{\ell_n} |\De_{i + j_m}^n X|^p = \sum_{i\De_n < s \le (i+k_n)\De_n} f(\De X(q)_s)  \Ind_B(s) \\ &~~~~~~~~~~~~~~~~~+\sum_{m=1}^{\ell_n} (A(n,q,i+j_m-1)_{(i+j_m)\De_n} + M(n,q,i+j_m-1)_{(i+j_m)\De_n})
\end{align*}
with $A(n,q,i+j_m-1)$ and $M(n,q,i+j_m-1)$ defined as above, but with $B$ being replaced by $((i+j_m-1)\De_n, (i+j_m)\De_n]$, also in the definition of $Y(q)$. Thus, as the respective sums over the jumps $f(\De X(q)_s)$ cancel, $A_n(q)$ becomes 
\begin{align} \label{410dar}
\nonumber \frac {2n}{\ell_n^2} &\sum_{i=0}^{n-k_n} \frac{1}{\binom{k_n}{\ell_n}} \sum_{1 \le j_1 < \ldots < j_{\ell_n} \le k_n} \Big(A(n,q,B_{i,j_1, \ldots, j_{\ell_n}}^n)_{(i+k_n)\De_n} + M(n,q,B_{i,j_1, \ldots, j_{\ell_n}}^n)_{(i+k_n)\De_n} \\&~~-\sum_{m=1}^{\ell_n} (A(n,q,i+j_m-1)_{(i+j_m)\De_n} + M(n,q,i+j_m-1)_{(i+j_m)\De_n})  \Big)^2 
\end{align}
and in order to show (\ref{410}) it becomes important to bound quantities like
\[
\E[(A(n,q,B)_{(i+k_n)\De_n})^2] \quad \text{and} \quad \E[(M(n,q,B)_{(i+k_n)\De_n})^2]  = \E[A'(n,q,B)_{(i+k_n)\De_n}].
\]
A Taylor expansion gives $|k(x,y)| \le C \left(|x| |y|^{p-1} + |y| |x|^{p-1}\right)$ as well as $|g(x,y)| \le  C \left(|x| |y|^{p-1} + y^2 |x|^{p-2}\right)$. From the boundedness conditions and integrability of $\ga(z)$ we obtain
\begin{align*}
|a(n,q,B)_u| &\le C \Ind_B(u) \left(q |Y(q)_{u-}|^{p-1} + |Y(q)_{u-}|^{p-2} + \al_q |Y(q)_{u-}| \right) \\ 
a'(n,q,B)_u	&\le C \Ind_B(u) \left(|Y(q)_{u-}|^{2p-2} + \al_q |Y(q)_{u-}|^2 \right) 
\end{align*}
for some sequence $\al_q$ with $\al_q \to 0$ as $q \to \infty$, and (15.2.22) in \cite{JacPro12} gives
\begin{align*}
\E\Big[\sup_{u \le (i+k_n) \De_n} |Y(q)_{u-}|^r \Big] \le C \left(q^r (\ell_n \De_n)^{r}  + (\ell_n \De_n)^{r/2} + \al_q (\ell_n \De_n)^{1 \wedge (r/2)}  \right)
\end{align*}
where we have used $|B| = \ell_n \De_n$. Therefore, 
\begin{align*}
\E\Big[\sup_{u \le (i+k_n) \De_n} a(n,q,B)_u^2 \Big]  \le C_q \ell_n \De_n
\end{align*}
and
\begin{align*}
\E\Big[\sup_{u \le (i+k_n) \De_n} a'(n,q,B)_u \Big] \le C_q (\ell_n \De_n)^2 + \al_q \ell_n \De_n.
\end{align*}
To summarize, 
\begin{align} \label{410a}
\E[(A(n,q,B)_{(i+k_n)\De_n})^2] \le (\ell_n \De_n)^2 \E\Big[\sup_{u \le (i+k_n) \De_n} a(n,q,B)_u^2 \Big] \le C_q(\ell_n \De_n)^3
\end{align}
and
\begin{align} \label{410b}
\nonumber \E[(M(n,q,B)_{(i+k_n)\De_n})^2] &= \E[A'(n,q,B)_{(i+k_n)\De_n}] \le \ell_n \De_n \E\Big[\sup_{u \le (i+k_n) \De_n} a'(n,q,B)_u \Big] \\ &\le  C_q (\ell_n \De_n)^3 + \al_q (\ell_n \De_n)^2.
\end{align}
Similar inequalities hold for $A(n,q,i+j_m-1)$ and $M(n,q,i+j_m-1)$, but with $\ell_n =1$. Then 
\begin{align} \label{410c}
\E\Big[ \big(\sum_{m=1}^{\ell_n} A(n,q,i+j_m-1)_{(i+j_m)\De_n}\big)^2 \Big] \le  C_q \ell_n^2 \De_n^3
\end{align}
and
\begin{align} \label{410d}
& \nonumber \E\Big[ \big(\sum_{m=1}^{\ell_n} M(n,q,i+j_m-1)_{(i+j_m)\De_n}\big)^2 \Big] \\ =& \sum_{m=1}^{\ell_n} \E\left[A'(n,q,i+j_m-1)_{(i+j_m)\De_n}\right] \le C_q \ell_n \De_n^3 + \al_q \ell_n \De_n^2.
\end{align}
From (\ref{410dar}) and the bounds in (\ref{410a})--(\ref{410d})  we obtain
\begin{align*}
\E\left[A_n(q) \right]  \le C \left(C_q \ell_n \De_n + \al_q\right),
\end{align*}
and the right hand side goes to zero as first $n \to \infty$ and then $q \to \infty$. This finishes the proof of (\ref{410}). 
 
The proof of (\ref{stabneg}) is complete by showing
\begin{align} \label{411}
\lim_{q \to \infty} \limsup_{n \to \infty} \P \left( B_n(q) > \eta \right) = 0
\end{align} 
which we will do under the assumption that (\ref{convnqfest}) holds. The proof of the latter claim will finish the entire section. Thus, let $\ka > 0$ be arbitrary. Then there exists $K > 0$ such that $\P(V \ge K) \le \ka$, and from the Portmanteau theorem we deduce 
\[
\limsup_{q \to \infty} \limsup_{n \to \infty} \P(V_n(q) \ge K) \le \limsup_{q \to \infty} \P(V(q) \ge K) \le \P(V \ge K) \le \ka.
\]
Let $\eps \le \frac{\eta}{3K}$. Then, using (\ref{ineqstd}), we obtain
\begin{align*}
&\P \left( B_n(q) > \eta \right) \le \P(\eps V_n(q) \Ind_{\{V_n(q) \ge K\}} > \eta/3) + \P(\eps V_n(q) \Ind_{\{V_n(q) < K\}} > \eta/3) \\ &+ \P \bigg(C_\eps \frac {n}{\ell_n(\ell_n-1)} \sum_{r \in P_q} \sum_{\al=1}^{k_n} \frac{1}{\binom{k_n}{\ell_n}} \sum_{\substack{1 \le j_1 < \ldots < j_{\ell_n} \le k_n \\ \{j_1, \ldots, j_{\ell_n}\} \cap \{ \al \} \neq \emptyset }}\big(|\De X_{S_r}|^p -\sum_{m = 1}^{\ell_n} |\De_{i_r - \al + j_m}^n X|^p \big)^2 > \eta/3 \bigg).
\end{align*}
For the first summand we have 
\[
\limsup_{q \to \infty} \limsup_{n \to \infty} \P(\eps V_n(q) \Ind_{\{V_n(q) \ge K\}} > \eta/3) \le  \P(V_n(q) \ge K) \le \ka,
\]
while for the second term 
\[
\P(\eps V_n(q) \Ind_{\{V_n(q) < K\}} > \eta/3) \le \P(\eps K > \eta/3) =0
\]
by construction. As $\ka$ was arbitrary (\ref{411}) follows, using $\ell_n \le 2(\ell_n-1)$ again, once we have shown 
\[
\lim_{q \to \infty} \limsup_{n \to \infty} \P \bigg( \frac {2n}{\ell_n^2} \sum_{r \in P_q} \sum_{\al=1}^{k_n} \frac{1}{\binom{k_n}{\ell_n}} \sum_{\substack{1 \le j_1 < \ldots < j_{\ell_n} \le k_n \\ \{j_1, \ldots, j_{\ell_n}\} \cap \{ \al \} \neq \emptyset }}\big(|\De X_{S_r}|^p -\sum_{m = 1}^{\ell_n} |\De_{i_r - \al + j_m}^n X|^p \big)^2 >\delta \bigg) = 0
\]
for any $\de > 0$, and we may assume to live on $\Om(n,q)$  without loss of generality. On this set  the decomposition 
\begin{align*}
& |\De_{i_r}^n X|^p - \sum_{m = 1}^{\ell_n} |\De_{i_r - \al + j_m}^n X|^p  \\ =& \bigg(\Big|\sum_{m=1}^{\ell_n} \De_{i_r-\al+j_m}^n X(q)\Big|^p -\sum_{m = 1}^{\ell_n} |\De_{i_r - \al + j_m}^n X(q)|^p \bigg) \\&-\Big|\sum_{m=1}^{\ell_n} \De_{i_r-\al+j_m}^n X(q)\Big|^p   \\&+ |\De_{i_r}^n X(q)|^p  \\ =& I(n,q,i_r,\al,j_1, \ldots, j_{\ell_n}) - II(n,	q,i_r,\al,j_1, \ldots, j_{\ell_n}) + III(n,q,i_r) 
\end{align*}
holds, because each interval $[i \De_n, (i+k_n)\De_n]$, $i=0, \ldots, n-k_n$, contains at most one jump of $N(q)$. We will now prove
\begin{align}
\label{bstep21} &\lim_{q \to \infty} \limsup_{n \to \infty} \P \Big( \frac {2n}{\ell_n^2} \sum_{r \in P_q} \sum_{\al=1}^{k_n} \frac{1}{\binom{k_n}{\ell_n}} \sum_{\substack{1 \le j_1 < \ldots < j_{\ell_n} \le k_n \\ \{j_1, \ldots, j_{\ell_n}\} \cap \{ \al \} \neq \emptyset }} I(n,q,i_r,\al,j_1, \ldots, j_{\ell_n})^2 > \de \Big) = 0, \\
\label{bstep22} &\lim_{q \to \infty} \limsup_{n \to \infty} \P \Big(\frac {2n}{\ell_n^2} \sum_{r \in P_q} \sum_{\al=1}^{k_n} \frac{1}{\binom{k_n}{\ell_n}} \sum_{\substack{1 \le j_1 < \ldots < j_{\ell_n} \le k_n \\ \{j_1, \ldots, j_{\ell_n}\} \cap \{ \al \} \neq \emptyset }} II(n,q,i_r,\al,j_1, \ldots, j_{\ell_n})^2 > \de \Big) = 0, \\
\label{bstep23} &\lim_{q \to \infty} \limsup_{n \to \infty} \P \Big(\frac{2n}{\ell_n}\sum_{r \in P_q}  III(n,q,i_r)^2 > \de \Big) = 0,
\end{align}
and again restricted to $\Om(n,q)$ if necessary. Note that the simplification in (\ref{bstep23}) is due to 
\begin{align} \label{binkoef}
\frac {n}{\ell_n^2} k_n \frac{1}{\binom{k_n}{\ell_n}} \binom{k_n-1}{\ell_n-1} = \frac{n}{\ell_n}.
\end{align}
Clearly, (\ref{bstep21}) is a simple consequence of (\ref{410}), and the proof of (\ref{bstep23}) is essentially the same as for (\ref{bstep22}), but with $\ell_n = 1$.

Thus, we will only prove (\ref{bstep22}), and we further introduce an auxiliary parameter $L \in \N$ and formally prove the equivalent convergence of (\ref{bstep22}) as first $n \to \infty$, then $L \to \infty$ and finally $q\to \infty$. Introducing the events $\Ind_{\{ |P_q| \le L \}}$ and $\Ind_{\{ |P_q| > L \}}$, where $|A|$ denotes the cardinality of a discrete set $A$, and from the fact that 
\[
\lim_{L \to \infty} \P(|P_q| > L) = 0
\]
for any fixed $q$, it is clear that (\ref{bstep22}) follows from 
\begin{align} \label{hilffix}
\lim_{n \to \infty} \P \Big(\frac {2n}{\ell_n^2} \sum_{\substack{r \in P_q\\|P_q| \le L}} \sum_{\al=1}^{k_n} \frac{1}{\binom{k_n}{\ell_n}} \sum_{\substack{1 \le j_1 < \ldots < j_{\ell_n} \le k_n \\ \{j_1, \ldots, j_{\ell_n}\} \cap \{ \al \} \neq \emptyset }} II(n,q,i_r,\al,j_1, \ldots, j_{\ell_n})^2 > \de \Big) = 0
\end{align}
for any fixed $q$ and $L$. As the sum over $r$ is then finite we may focus on a single arbitrary index $i_r$, and by properties of a Poisson measure we can also drop the dependence on the jumps of $\Ind_{L_q} \star \mu$ and simply write $i$. With the notation (\ref{defyq}) we have 
\[
II(n,q,i,j_1, \ldots, j_{\ell_n}) \le \left| \int_{(i-k_n) \De_n}^{(i+k_n)\De_n} \Ind_{B_{i-\al,j_1, \ldots, j_{\ell_n}}^n}(s) dX(q)_s \right|^{p}.
\]
By definition $X(q)$ consists of three terms, and we will discuss each of them separately. The first two are easier to deal with, and we have
\[
\E \bigg[ \Big| \int_{(i-k_n) \De_n}^{(i+k_n)\De_n} \Ind_{B_{i-\al,j_1, \ldots, j_{\ell_n}}^n}(s) b(q)_s ds \Big|^{2p} \bigg] \le C_q (\ell_n \De_n)^{2p}
\]
and
\begin{align} \label{contmartp}
\E \bigg[ \Big| \int_{(i-k_n) \De_n}^{(i+k_n)\De_n} \Ind_{B_{i-\al,j_1, \ldots, j_{\ell_n}}^n}(s) \si_s dW_s \Big|^{2p} \bigg] \le C (\ell_n \De_n)^{p}.
\end{align}
Together with (\ref{binkoef}) it is clear that (\ref{hilffix}) follows from  
\[
\E\Big[\Big(\frac {2n}{\ell_n^2} \sum_{\al=1}^{k_n} \frac{1}{\binom{k_n}{\ell_n}} \sum_{\substack{1 \le j_1 < \ldots < j_{\ell_n} \le k_n \\ \{j_1, \ldots, j_{\ell_n}\} \cap \{ \al \} \neq \emptyset }} \Big(\int_{(i-k_n) \De_n}^{(i+k_n)\De_n}  \int \delta(s,z) \Ind_{B}(s) \mathds{1}_{\{\gamma(z)\leq 1/q\}}(\mu-\nu)(ds,dz) \Big)^{2p} \Big) \wedge 1 \Big] \to 0
\]
where we again use the notation $B=B_{i-\al,j_1, \ldots, j_{\ell_n}}^n$. For any $0 < \eps < 1$ and any $t \ge (i-k_n) \De_n$ we decompose the above integral into three terms and set
\begin{align*}
\widetilde N(\eps)_t &= \int_{(i-k_n) \De_n}^{t} \int \Ind_{\{\ga(z) > \eps\}} \mu(ds,dz), \\
\widetilde M(\eps)_t &= \int_{(i-k_n) \De_n}^{t} \int \Ind_{\{\ga(z) \le \eps\}}  \delta(s,z) \Ind_{B}(s) \mathds{1}_{\{\gamma(z)\leq 1/q\}}(\mu-\nu)(ds,dz), \\
\widetilde B(\eps)_t &= -\int_{(i-k_n) \De_n}^{t} \int \Ind_{\{\ga(z) > \eps}\}  \delta(s,z) \Ind_{B}(s) \mathds{1}_{\{\gamma(z)\leq 1/q\}} \la(dz) ds.
\end{align*}
By integrability of $\ga^2$ we have
\begin{align*}
\P(\widetilde N(\eps)_{(i+k_n)\De_n} \ge 1) \le \E[\widetilde N(\eps)_{(i+k_n)\De_n}] = \E \left[\int_{(i-k_n) \De_n}^{(i+k_n)\De_n} \int \Ind_{\{\ga(z) > \eps\}} \la(dz) ds  \right] \le C \frac{k_n \De_n}{\eps^2},
\end{align*}
and similarly we can deduce $|\widetilde B(\eps)_{(i+k_n)\De_n}| \le C \frac{\ell_n \De_n}{\eps}.$ Finally, from Lemma 2.1.5 in \cite{JacPro12} we obtain 
\begin{align} \label{ineqm}
&\E[|\widetilde M(\eps)_{(i+k_n)\De_n}|^{2p}] \\ \le& C \bigg(\E \Big[\int_{(i-k_n) \De_n}^{(i+k_n)\De_n} \int \Ind_{\{\ga(z) \le \eps\}}   |\delta(s,z)|^{2p} \Ind_{B}(s) \mathds{1}_{\{\gamma(z)\leq 1/q\}} \la(dz) ds \Big] \nonumber \\ \nonumber &~~~~+ \E \Big[\Big(\int_{(i-k_n) \De_n}^{(i+k_n)\De_n} \int \Ind_{\{\ga(z) \le \eps\}}   |\delta(s,z)|^{2} \Ind_{B}(s) \mathds{1}_{\{\gamma(z)\leq 1/q\}} \la(dz) ds \Big)^p \Big] \bigg)
\\ \le& C \bigg(\eps^{2p-2} \ell_n \De_n \int \ga(z)^{2} \la(dz)  + (\ell_n \De_n)^p \Big(\int \ga(z)^{2} \la(dz)\Big)^p\bigg) \le C (\eps^{2p-2} \ell_n \De_n   + (\ell_n \De_n)^p). \nonumber
\end{align}
Then 
\begin{align*} 
& \E\Big[\Big(\frac {2n}{\ell_n^2} \sum_{\al=1}^{k_n} \frac{1}{\binom{k_n}{\ell_n}} \sum_{\substack{1 \le j_1 < \ldots < j_{\ell_n} \le k_n \\ \{j_1, \ldots, j_{\ell_n}\} \cap \{ \al \} \neq \emptyset }} \Big(\int_{(i-k_n) \De_n}^{(i+k_n)\De_n}  \int \delta(s,z) \Ind_{B}(s) \mathds{1}_{\{\gamma(z)\leq 1/q\}}(\mu-\nu)(ds,dz) \Big)^{2p} \Big) \wedge 1 \Big]
\\ &\le \P(\widetilde N(\eps)_{(i+k_n)\De_n} \ge 1) + C_p \frac {n}{\ell_n^2} \sum_{\al=1}^{k_n} \frac{1}{\binom{k_n}{\ell_n}} \sum_{\substack{1 \le j_1 < \ldots < j_{\ell_n} \le k_n \\ \{j_1, \ldots, j_{\ell_n}\} \cap \{ \al \} \neq \emptyset }} \E[|\widetilde B(\eps)_{(i+k_n)\De_n}|^{2p} + |\widetilde M(\eps)_{(i+k_n)\De_n}|^{2p}]
\\ &\le C_p \left( \frac{k_n \De_n}{\eps^2} + \left({\ell_n \De_n}\right)^{2p-1} \eps^{-2p} + \eps^{2p-2} + (\ell_n \De_n)^{p-1}\right),
\end{align*}
where we have used (\ref{binkoef}). Choosing $\eps_n \to 0$ small enough then ends the proof of (\ref{hilffix}).


We will finish the proof by showing (\ref{convnqfest}), for which we use the following Taylor expansion for $f(x) = |x|^p$: On $\Om(n,q)$ we have 
\begin{align*}
&\Big|\sum_{m=1}^{\ell_n} \De_{i_r-\al+j_m}^n X \Big|^p - |\De_{i_r}^n X|^p = f'(\De X_{S_r})  \sum_{\substack{m = 1 \\ j_m \neq \al}}^{\ell_n} \De_{i_r-\al+j_m}^n X(q) \\ &~~~~+ f''(\ka_{i_r}^n) \De_{i_r}^n X(q)  \sum_{\substack{m = 1 \\ j_m \neq \al}}^{\ell_n} \De_{i_r-\al+j_m}^n X(q) + \frac 12 f''(\xi_{i_r - \al, j_1, \ldots, j_{\ell_n}}^n) \Big| \sum_{\substack{m = 1 \\ j_m \neq \al}}^{\ell_n} \De_{i_r-\al+j_m}^n X(q)\Big|^2
\end{align*}
for some intermediate $\ka_{i_r}^n$ between $\De_{i_r}^n X$ and $\De X_{S_r}$ and $\xi_{i_r - \al, j_1, \ldots, j_{\ell_n}}^n$ between $\sum_{m=1}^{\ell_n} \De_{i_r-\al+j_m}^n X $ and $\De_{i_r}^n X$. On $\Om(n,q)$ both are bounded by $C_q$. Obviously, one can show 
\[
\frac {2n}{\ell_n^2} \sum_{r \in P_q} \sum_{\al=1}^{k_n} \frac{1}{\binom{k_n}{\ell_n}} \sum_{\substack{1 \le j_1 < \ldots < j_{\ell_n} \le k_n \\ \{j_1, \ldots, j_{\ell_n}\} \cap \{ \al \} \neq \emptyset }} \Big| \sum_{\substack{m = 1 \\ j_m \neq \al}}^{\ell_n} \De_{i_r-\al+j_m}^n X(q)\Big|^4 \pn 0
\] 
as $n \to \infty$ for any fixed $q$ along the same lines as the ones from the proof of (\ref{bstep22}) with $p=2$, and similarly 
\[
\frac {2n}{\ell_n^2} \sum_{r \in P_q} \sum_{\al=1}^{k_n} \frac{1}{\binom{k_n}{\ell_n}} \sum_{\substack{1 \le j_1 < \ldots < j_{\ell_n} \le k_n \\ \{j_1, \ldots, j_{\ell_n}\} \cap \{ \al \} \neq \emptyset }} |\De_{i_r}^n X(q)|^2 \Big| \sum_{\substack{m = 1 \\ j_m \neq \al}}^{\ell_n} \De_{i_r-\al+j_m}^n X(q)\Big|^2 \pn 0.
\] 
Lemma \ref{lemyoung} then suggests that we only need to prove 
\[
\frac {n}{\ell_n(\ell_n-1)} \sum_{r \in P_q} \sum_{\al=1}^{k_n} \wY_{r,\al}^{(n)} \pn V(q)
\]
where
\[
\wY_{r,\al}^{(n)} = \frac{1}{\binom{k_n}{\ell_n}} \sum_{\substack{1 \le j_1 < \ldots < j_{\ell_n} \le k_n \\ \{j_1, \ldots, j_{\ell_n}\} \cap \{ \al \} \neq \emptyset }} f'(\De X_{S_r})^2 \Big| \sum_{\substack{m = 1 \\ j_m \neq \al}}^{\ell_n} \De_{i_r-\al+j_m}^n X(q)\Big|^2.
\]
The penultimate step is yet another application of Lemma \ref{lemyoung}, namely to first prove 
\[
\frac {2n}{\ell_n^2} \sum_{r \in P_q} \sum_{\al=1}^{k_n} \frac{1}{\binom{k_n}{\ell_n}} \sum_{\substack{1 \le j_1 < \ldots < j_{\ell_n} \le k_n \\ \{j_1, \ldots, j_{\ell_n}\} \cap \{ \al \} \neq \emptyset }} \Big( \sum_{\substack{m = 1 \\ j_m \neq \al}}^{\ell_n} (\De_{i_r - \al + j_m}^n X(q) - \si_{S_r} \De_{i_r - \al + j_m}^n W)\Big)^2 \pn 0
\] 
as $n \to \infty$ for any fixed $q$ and to use boundedness of the jumps of $N(q)$ by some $C_q$. This proof also works in the same way as (\ref{bstep22}) with $p=1$, but with two differences: First, instead of \eqref{contmartp} we discuss
\begin{align*}
\E \bigg[ \Big| \int_{(i_r-k_n) \De_n}^{(i_r+k_n)\De_n} \Ind_{B_{i-\al,j_1, \ldots, j_{\ell_n}}^n}(s) (\si_s -  \si_{S_r}) dW_s \Big|^{2} \bigg] \le \ell_n \De_n \E[\sup_{|u| \le k_n \De_n} |\si_{S_r - u} - \si_{S_r}|^2],
\end{align*}
and we apply additionally continuity of $\si$ in $S_r$ plus dominated convergence, and second the upper bound in (\ref{ineqm}) now becomes $\ell_n \De_n \int\Ind_{\{\ga(z) \le \eps\}}  \ga(z)^{2} \la(dz)$. Therefore, from Lemma \ref{lemyoung} it is sufficient to prove convergence in probability of 
\begin{align} \label{stepult}
\frac {n}{\ell_n (\ell_n-1)} \sum_{r \in P_q} (f'(\De X_{S_r}))^2 \si^2_{S_r} \sum_{\al=1}^{k_n} \frac{1}{\binom{k_n}{\ell_n}} \sum_{\substack{1 \le j_1 < \ldots < j_{\ell_n} \le k_n \\ \{j_1, \ldots, j_{\ell_n}\} \cap \{ \al \} \neq \emptyset }} \Big|  \sum_{\substack{m = 1 \\ j_m \neq \al}}^{\ell_n} \De_{i_r-\al+j_m}^n W \Big|^2
\end{align}
to $V(q)$ as $n \to \infty$. Using $f'(x) = p x^{p-1}$ and \eqref{binkoef2} we are left to show
$\frac {1}{k_n}  \sum_{\al=1}^{k_n} Z_{i,\al}^n \pn 0$ for any fixed $i$, where 
\[
Z_{i,\al}^n= \frac{1}{\binom{k_n-1}{\ell_n-1}} \sum_{\substack{1 \le j_1 < \ldots < j_{\ell_n-1} \le k_n \\ \{j_1, \ldots, j_{\ell_n-1}\} \cap \{ \al \} = \emptyset }} \bigg( \frac{n}{\ell_n-1}\Big|\sum_{m = 1}^{\ell_n-1} \De_{i - \al + j_m}^n W \Big|^2 - 1 \bigg).
\]
Note that we can again drop the dependence on $r$ by properties of a Poisson random measure. Using $\E[Z_{i_r,\al}^n] = 0$ and 
\[
\Var \Big(\frac {1}{k_n}  \sum_{\al=1}^{k_n} Z_{i,\al}^n \Big) = \frac 1{k_n^2} \sum_{\al_1, \al_2=1}^{k_n} \Cov(Z_{i,\al_1}^n, Z_{i,\al_2}^n) \le \Big(\frac {1}{k_n}  \sum_{\al=1}^{k_n} \sqrt{ \Var(Z_{i,\al}^n) } \Big)^2
\]
we are left to show $\Var(Z_{i,\al}^n) \le \eta_n \to 0$. In distribution, $Z_{i,\al}^n$ equals the U statistic 
\[
U_n = \frac{1}{\binom{k_n-1}{\ell_n-1}} \sum_{1 \le j_1 < \ldots < j_{\ell_n-1} \le k_n-1} \bigg(\Big|\sum_{m = 1}^{\ell_n-1} N_{j_m}^n \Big|^2 - 1 \bigg)
\]
for i.i.d.\ standard normal $N_i$. Using Theorem 1.2.3 in \cite{denker1985} again we obtain 
\[
\Var(U_n) \le C \frac{\ell_n-1}{k_n-1} \to 0
\]
which finishes the proof for $V_n$.

We will finally discuss the necessary changes for $\wV_n$ and $\WV_n$, and this time the entire proof goes through in exactly the same way when $\ell_n = 2$. For $\ell_n = k_n$ the proof of (\ref{410}) goes through without any changes, whereas for (\ref{411}) we cannot apply (\ref{convnqfest}) because we do not have convergence in probability in the end. Nevertheless, we only use (\ref{convnqfest}) in an application of the Portmanteau theorem, and this goes through under weak convergence as well. So we only need to discuss the stable convergence of (\ref{convnqfest}), as (\ref{convq}) finally follows from monotone convergence again.

The proof of (\ref{convnqfest}) can always be reproduced until one arrives at (\ref{stepult}) which, because of $k_n \to \infty$, becomes
\begin{align*} 
& \sum_{r \in P_q} (f'(\De X_{S_r}))^2 \si^2_{S_r} \frac 1{k_n} \sum_{\al=1}^{k_n} \frac{n}{k_n} \Big|  \sum_{m = 1}^{k_n} \De_{i_r-\al+m}^n W \Big|^2 \big(1+ o_\P(1)\big) \\ =& \sum_{r \in P_q} (f'(\De X_{S_r}))^2 \si^2_{S_r} w_{n,r} \big(1+ o_\P(1)\big)
\end{align*} 
with
\[
w_{n,r} = \frac{n}{k_n^2}  \sum_{j=0}^{k_n-1} (W_{(i_r+k_n-j)\De_n} - W_{(i_r-j) \De_n})^2. 
\]
The final step therefore is to prove the stable convergence
\[
(w_{n,r})_{r \ge 1} \tols (1+R_r)_{r \ge 1},
\]
which follows as in the proof of Theorem 4.3.1 in \cite{JacPro12} and can be traced back to convergence in distribution of each fixed $w_{n,r}$ to $1+R_r$. This latter convergence is granted using Theorem 1 in \cite{wushao2007}. Note that this result is concerned with convergence in distribution to a limiting normal distribution. Note, however, that their condition (15) is not satisfied in our situation. Nevertheless, convergence in distribution still holds, see the comment following their Theorem 1, but the limiting distribution remains unspecified. \qed

\bibliographystyle{chicago}
\bibliography{bibliography}
\end{document}